\documentclass[a4paper]{amsart}\usepackage{amsfonts,amssymb,bbm,graphicx}
\usepackage{enumerate}

\def\C{\Bbb{C}}\def\k{\mathbbm{k}}
\def\bI{{\bar{I}}}\def\N{\Bbb{N}}
\def\R{\Bbb{R}}

\def\kpd{{$\k$-polynomially-defined}}

\def\bl{\langle}\def\br{\rangle}\def\suml{\sum\limits}\def\oplusl{\mathop\oplus\limits}\def\liml{\lim\limits}

\def\capl{\mathop\cap\limits}\def\bigcapl{\mathop\bigcap\limits}
\def\prodl{\mathop\prod\limits}

\newcommand{\quotient}[2]{{\left.\raisebox{1.6ex}{$#1$}\!\!\!\!\!{\scalebox{2}{\ensuremath\diagup}}
\!\!\!\!\!\raisebox{-1ex}{$#2$}\right.}}
\newcommand{\quots}[2]{{\footnotesize\left.\raisebox{0.9ex}{$#1$}\!\!\! \diagup \!\!\!\raisebox{-0.9ex}{$#2$}\right.}}

\def\tA{\widetilde{A}}
\def\tJ{\tilde{J}}
\def\tU{\tilde{U}}

\def\hA{{\widehat{A}}}\def\hF{\hat{F}}
\def\hi{{\widehat{i}}}
\def\hM{{\hat{M}}}\def\hR{{\widehat{R}}}\def\hu{\hat{u}}\def\hv{\hat{v}}

\def\de{\delta}\def\De{\Delta}

\def\la{\lambda}\def\La{\Lambda}\def\Si{\Sigma}

\def\cA{\mathcal A}\def\ca{\mathfrak a}

\def\cP{\mathcal P}

\def\cm{{\frak m}} 

\def\ux{\underline{x}}

\def\one{{1\hspace{-0.1cm}\rm I}}\def\zero{\mathbb{O}}

\newcommand{\ber}{\begin{array}{l}}\newcommand{\eer}{\end{array}}
\newcommand{\bpm}{\begin{pmatrix}}\newcommand{\epm}{\end{pmatrix}}
\newcommand{\bM}{\begin{matrix}}\newcommand{\eM}{\end{matrix}}
\newcommand{\bee}{\begin{enumerate}}\newcommand{\eee}{\end{enumerate}}
\newcommand{\bei}{\begin{itemize}}\newcommand{\eei}{\end{itemize}}

\def\wrt{with respect to }\def\sset{\subset}\def\sseteq{\subseteq}\def\ssetneq{\subsetneq}\def\smin{\setminus}

\def\Mat{Mat_{m\times n}(R)}\def\Matm{Mat_{m\times n}(\cm)}

\def\bull{\vrule height .9ex width .9ex depth -.1ex }

\newtheorem{Lemma}{Lemma}[section]\newcommand{\bel}{\begin{Lemma}}\newcommand{\eel}{\end{Lemma}}
\newtheorem{Theorem}[Lemma]{Theorem}\newcommand{\bthe}{\begin{Theorem}}\newcommand{\ethe}{\end{Theorem}}
\newtheorem{Proposition}[Lemma]{Proposition}\newcommand{\bprop}{\begin{Proposition}}\newcommand{\eprop}{\end{Proposition}}
\newtheorem{Corollary}[Lemma]{Corollary}\newcommand{\bcor}{\begin{Corollary}}\newcommand{\ecor}{\end{Corollary}}
\newtheorem{Definition}[Lemma]{Definition}\newcommand{\bed}{\begin{Definition}}\newcommand{\eed}{\end{Definition}}
\newtheorem{Definition-Proposition}[Lemma]{Definition-Proposition}

\def\bpr{~\\{\em Proof.\ }}\newcommand{\epr}{$\bull$\\}
\newtheorem{Remark}[Lemma]{Remark}\newcommand{\beR}{\begin{Remark}\rm}\newcommand{\eeR}{\end{Remark}}
\newtheorem{Example}[Lemma]{Example}\newcommand{\bex}{\begin{Example}\rm}\newcommand{\eex}{\end{Example}}
\newtheorem{Property}[Lemma]{Property}\newcommand{\bproperty}{\begin{Property}\rm}\newcommand{\eproperty}{\end{Property}}
\newtheorem{Problem}[Lemma]{Problem}\newcommand{\bprob}{\begin{Problem}\rm}\newcommand{\eprob}{\end{Problem}}

\newcommand{\bet}{\begin{tabular}{cccccccc}}\newcommand{\eet}{\end{tabular}}
\newcommand{\beq}{\begin{equation}}\newcommand{\eeq}{\end{equation}}

\newcommand{\bin}[2]{\binom{#1}{#2}}
\newcommand{\cornA}[1]{\underline{A\ }\!|\!\!|_{#1}}

\newcommand\isom{\xrightarrow{\,\smash{\raisebox{-0.65ex}{\ensuremath{\scriptstyle\sim}}}\,}}

\title[]{F\MakeLowercase{inite determinacy of matrices over local rings.}
\\T\MakeLowercase{angent modules to the miniversal deformation for} $R$-\MakeLowercase{linear group actions.}}
\author{G\MakeLowercase{enrich} B\MakeLowercase{elitskii and} D\MakeLowercase{mitry} K\MakeLowercase{erner}}
\address{Department of Mathematics, Ben Gurion University of the Negev, P.O.B. 653, Be'er Sheva 84105, Israel.}
\email{genrich@math.bgu.ac.il}
\email{dmitry.kerner@gmail.com}
\date{\today}
\thanks{D.K. was supported by the  Israel Science Foundation (grant No. 844/14).}

\begin{document}\maketitle
\begin{abstract}
We consider matrices with entries in a local ring, $\Mat$. Fix a group action, $G\circlearrowright\Mat$, and a subset of allowed deformations, $\Si\sseteq\Mat$.
 The standard question in Singularity Theory is the finite-$(\Si,G)$-determinacy of matrices. Finite determinacy implies algebraizability and
  is equivalent to a stronger notion: stable algebraizability.

 In our previous work this determinacy question was reduced to the study of the tangent spaces $T_{(\Si,A)}$, $T_{(GA,A)}$, and their quotient, the tangent module to the miniversal deformation, $T^1_{(\Si,G,A)}=\quots{T_{(\Si,A)}}{T_{(GA,A)}}$.
 In particular, the order of determinacy is controlled by the annihilator of this tangent module, $ann(T^1_{(\Si,G,A)})$.

In this work we  study this tangent module for the group action $GL(m,R)\times GL(n,R)\circlearrowright\Mat$ and  various natural subgroups of it.

We obtain  ready-to-use criteria of determinacy for deformations of
 (embedded) modules, (skew-)symmetric forms, filtered modules, filtered morphisms of filtered
modules, chains of modules and others.
\end{abstract}\setcounter{secnumdepth}{6} \setcounter{tocdepth}{1}

\section{Introduction}

\subsection{Setup}\label{Sec.Intro.Setup}
Let $R$ be a (commutative, associative) local ring over some base field $\k$ of zero characteristic.
 (The simplest examples are $\k[[\ux]]$, $\C\{\ux\}$, and $C^\infty(\R^p,0)$.)
Denote by $\cm\sset R$ the maximal ideal and by
 $\Mat$ the $R$-module of $m\times n$ matrices with entries in $R$. We always assume $m\le n$, otherwise one
can transpose the matrix.

\

Various groups act on $\Mat$, we consider mostly the following actions.\vspace{-0.1cm}
\bex\label{Ex.Intro.Typical.Groups} Let $A\in \Mat$.
\bee[i.]
\item
 The left multiplications $G_l:=GL(m,R)$, the right and the two-sided multiplications $G_r:=GL(n,R)$,  $G_{lr}:=G_l\times G_r$. Matrices
  considered up to $G_r$-transformations correspond to the embedded modules,  $Im(A)\sset R^{\oplus m}$. Matrices
   considered up to $G_{lr}$-transformations correspond to the non-embedded modules,  $coker(A)=\quots{R^{\oplus m}}{Im(A)}$.
\item The congruence, $G_{congr}=GL(m,R)\circlearrowright Mat_{m\times m}(R)$, acts by $A\to UAU^t$. Matrices considered up
 to the congruence correspond to the bilinear/symmetric/skew-symmetric forms. The conjugation, $G_{conj}=GL(m,R)\circlearrowright Mat_{m\times m}(R)$,
  acts by $A\to UAU^{-1}$. This action is important in representation theory.
\item Fix some ordered compositions, $m=\suml^k_{i=1}m_i$, $n=\suml^k_{j=1}n_j$. Take the corresponding block-structure on the matrix, i.e. $A=\{A_{ij}\}_{1\le i,j\le k}$, $A_{ij}\in Mat_{m_i\times n_j}(R)$.
Denote by $Mat^{up}_{\{m_i\}\times\{n_j\}}(R)$ the set of the upper-block-triangular matrices, i.e. $A_{ij}=\zero$ for $i>j$.
 Thus $Mat^{up}_{\{m_i\}\times\{n_j\}}(R)$
  is a free direct summand of $\Mat$.
 Accordingly we consider the groups
\beq
 G^{up}_l:=G_l\cap Mat^{up}_{\{m_i\}\times\{m_j\}}(R),\quad
 G^{up}_r:=G_r\cap Mat^{up}_{\{n_i\}\times\{n_j\}}(R)
 \ \text{ and } \ G^{up}_{lr}:=G^{up}_l\times G^{up}_r.
\eeq
  The upper-block-triangular matrices considered up to the $G^{up}_r$-equivalence correspond to the submodules of $R^{\oplus m}$
   filtered by chains of direct summands of $R^{\oplus m}$. The upper-block-triangular matrices considered up to the $G^{up}_{lr}$-equivalence correspond to filtered
    homomorphisms of filtered free modules. (See \S \ref{Sec.Proofs.Upper.Triangular} for more detail.)
\eee
\eex

\

We study deformations of matrices. In applications one often deforms a matrix not inside the whole $\Mat$ but only inside a ``deformation subspace" (a subset of prescribed deformations), $A\rightsquigarrow A+B$, $A+B\in\Si\sseteq\Mat$.

\

In this paper the subset $\Si-\{A\}\sseteq\Mat$ is a submodule. Besides the trivial choice $\Si=\Mat$, we mostly work with the following submodules.
\vspace{-0.1cm} 
\bex\label{Ex.Intro.Typical.Def.Spaces}
\bee[i.]
\item The congruence of (skew-)symmetric matrices, $A\stackrel{G_{congr}}{\sim}UAU^t$, $U\in GL(m,R)$, preserves the (skew-) symmetry.
Thus it is natural to deform the matrix by the (skew-)symmetric matrices only, i.e. $\Si=Mat^{sym}_{m\times m}(R)$ or $\Si=Mat^{skew-sym}_{m\times m}(R)$.
\item  Similarly, it is natural to deform an upper-block-triangular matrix only inside $Mat^{up}_{\{m_i\}\times \{n_j\}}(R)$.
\item Sometimes one deforms only by ``higher order terms", e.g. $\Si=\{A\}+Mat_{m\times m}(J)$ or $\Si=\{A\}+Mat^{sym}_{m\times m}(J)$, for some ideal $J\sset R$.
\eee
\eex

\

 Recall that over a local ring any matrix is $G_{lr}$-equivalent to a block-diagonal, $A\stackrel{G_{lr}}{\sim}\one\oplus \tA$,
where  all the entries of $\tA$ lie in the maximal ideal $\cm$, i.e. vanish at the origin of $Spec(R)$.
 Similar statements hold for (skew-)symmetric matrices \wrt $G_{congr}$, see Proposition \ref{Thm.Background.Chip.off.Unit}.
 This splitting is natural in various senses and is standard in Commutative Algebra,  Singularity Theory and other fields.
 Often it is $\tA$ that carries all the essential information.
 Therefore we often assume $A|_0=\zero$, i.e. $A\in\Matm$.

\subsection{The tangent spaces}
Fix an action $G\circlearrowright\Mat$, a deformation space $\Si\sseteq\Mat$, and a matrix $A\in\Si$. We assume that the
 orbit $GA$ and the deformation space possess well defined tangent spaces at $A$, which are $R$-modules.
 While the general definitions/conditions are technical (see \cite{Belitski-Kerner}), in all the cases relevant to our paper
 the tangent spaces are ``the expected ones",  see \S\ref{Sec.Background.Tangent.Spaces}.

The standard approach of deformation theory is to establish the existence of the miniversal (semi-universal) deformation and, when the later exists, to understand/to compute its tangent cone.
Accordingly one passes from the study of the germs $(GA,A)\sseteq(\Si,A)$, to the study of the tangent spaces, $T_{(GA,A)}\sseteq T_{(\Si,A)}$.

Much of the information about the deformation problem is encoded in the quotient module
\beq
T^1_{(\Si,G,A)}:=\quots{T_{(\Si,A)}}{T_{(GA,A)}}.
\eeq
 This $R$-module is the tangent space to the miniversal deformation, when the later exists and is smooth. It often controls the deformation theory of modules, (skew-)symmetric forms and other objects.
In Singularity Theory such a module is known as the Tjurina algebra for the contact equivalence, and the Milnor algebra for the right equivalence, \cite{AGLV}.

\

This module $T^1_{(\Si,G,A)}$ is the main object of our study.

\subsection{Finite and infinite determinacy}
 The $(\Si,G)$-order of determinacy of $A$ is  the minimal number $k\le\infty$ satisfying: if $A,A_1\in\Si$ and $jet_k(A)=jet_k(A_1)$ then $A_1\in GA$. Here $jet_k$ is the projection $\Mat\stackrel{jet_k}{\to}Mat_{m\times n}(\quots{R}{\cm^{k+1}})$. More precisely:
\bed\label{Def.Intro.Determinacy.Classical}
$ord^\Si_G(A):=min\Big\{k|\ \Si\cap\big(\{A\}+Mat_{m\times n}(\cm^{k+1})\big)\sseteq GA\Big\}\le\infty.$
\eed
(We assume that the minimum is taken here over a non-empty set.)

\

If $ord^\Si_G(A)<\infty$ then $A$ is called {\em finitely-$(\Si,G)$-determined}.
 In words: $A$ is determined (up to $G$-equivalence) by a finite number of terms of its Taylor expansion at the origin.

Usually $R$  contains a polynomial subring, $\quots{\k[\ux]}{I}$, with the property $\cm=R\cdot(\ux)$, i.e. the elements of $R$ are power series/functions in $\ux$.
 Then an immediate consequence of finite determinacy is the {\em algebraization}:  $A$ is $G$-equivalent to a matrix of polynomials in $\ux$.
   Moreover, the order of determinacy gives an upper bound on the degrees of polynomials. We show that the finite determinacy is equivalent to
  the stable algebraizability (a strengthening of the ordinary algebraizability), see \S\ref{Sec.Results.Algebraization} and \S\ref{Sec.Background.Fin.Determin.vs.Stable.Algebraizability}.

 If $ord^\Si_G(A)=\infty$ the matrix is called {\em infinitely determined}.
 This condition is empty when $\cm^\infty:=\capl_{k>0}\cm^k=\{0\}$.
 But for some rings $\{0\}\neq \cm^\infty\sset R$, then even the infinite determinacy is a non-trivial  property.
  (The main example is the ring of germs of smooth functions, $R=C^{\infty}(\R^p,0)$, and its sub-quotients.) In this case $A$ is determined (up to $G$-equivalence) by its image under the $\cm$-adic completion, $\hA\in Mat_{m\times n}(\hR)$, i.e. its full Taylor expansion at the origin.

\

The finite and infinite determinacies are fundamental notions of Singularity Theory. For Algebraic Geometry/Commutative Algebra the finite determinacy
 means (roughly) ``the deformation theory is essentially finite dimensional''.
More generally, as the determinacy expresses the ``minimal stability", it is important in any area dealing with matrices over rings
(or matrix families or matrices depending on parameters).

 See \S\ref{Sec.Results.Remarks} for further motivation and \S\ref{Sec.Results.Relation.to.Singularity.Theory} for a brief overview and the
  relation of our work to the known results.

\

In \cite{Belitski-Kerner} we have reduced the study of determinacy to the understanding of the support/annihilator of the module
 $T^1_{(G,\Si,A)}$. The matrix $A$ is (in)finitely determined iff $\cm^N T_{(\Si,A)}\sseteq T_{(GA,A)}$ for some $N\le\infty$, equivalently
  $\cm^N T^1_{(G,\Si,A)}=\{0\}$ or $\cm^N\sseteq ann(T^1_{(G,\Si,A)})$. The order of determinacy is
  fixed by the annihilator $ann(T^1_{(G,\Si,A)})$, see  \S\ref{Sec.Background.BK1} for more detail.

\

For Artinian rings the finite determinacy is not interesting, e.g. if $\cm^{N+1}=\{0\}$ then any matrix
  over $R$ is trivially $N$-determined. Therefore
 in this paper we often assume that $R$ is not Artinian, i.e. $\cm^N\neq\{0\}$ for any $N<\infty$. Geometrically this means: $dim(Spec(R))>0$.

\subsection{Contents of the paper}
The results of \cite{Belitski-Kerner} reduce the study of (in)finite determinacy to the computation of the support of $T^1_{(\Si,G,A)}$, i.e.
 the annihilator ideal $ann(T^1_{(\Si,G,A)})$.
In the current paper we study the module $T^1_{(\Si,G,A)}$ and its support for the $R$-linear actions $G\circlearrowright \Mat$ of example
 \ref{Ex.Intro.Typical.Groups}.

\

One could begin from a completely general (set-theoretic) action $G\circlearrowright \Mat$.
 However, if an action $G\circlearrowright\Mat$ is $R$-linear (though does not necessarily preserve the row/column structure)
 and sends the degenerate matrices to degenerate (quite a reasonable assumption!) then $G\sseteq G_{lr}$. This property is proved in
 \S\ref{Sec.Background.Largest.Reasonable.Group}, it  belongs to the area of ``preserver problems on matrices", exhaustively studied for matrices over
  fields, but less known for matrices over rings.  Therefore example
 \ref{Ex.Intro.Typical.Groups} contains the main ``reasonable" cases.

\

 We compute (or at least bound) the annihilator $ann(T^1_{(\Si,G,A)})$, the main results are stated in \S\ref{Sec.Results}.
  This gives
 the ready-to-use criteria to determine (or at least to bound) $ord^\Si_G(A)$ for  various group-actions.
  We observe the standard dichotomy: for a given data $(m,n,dim(R),G\circlearrowright\Si\sseteq\Mat)$ either
   there are no finitely determined matrices in $\Matm\cap\Si$, or the generic finite determinacy holds.

The proofs of these results, corollaries/examples and further developments are in \S\ref{Sec.Proofs}.
In many cases we strengthen (generalize and quantify) both the classical and the relatively recent  results.


\

In \S\ref{Sec.Background} we collect the needed background material: determinantal ideals and Pfaffian ideals, annihilator-of-cokernel,
 the integral closure of ideals and modules, the matrix preservers, tangent spaces to the group orbits, approximation properties and the
 main result of \cite{Belitski-Kerner}.

\subsection{Acknowledgements} Many thanks to A.Fernandez-Boix, R.-O. Buchweitz, V. Grandjean, G.-M. Greuel, V. Kodiyalam and M. Leyenson for
helpful/informative discussions and to L. Moln\'{a}r, V.V.
Sergeichuk for the much needed references on the ``matrix
preservers" problems. Finally, we thank the anonymous referee for careful reading of the paper and valuable suggestions.

\section{The main results}\label{Sec.Results}
Below we state the main theorems and corollaries. Their proofs and examples are in \S\ref{Sec.Proofs}.

Unless stated otherwise, $R$ is just a commutative unital ring, not necessarily local or Noetherian.

\subsection{Notations and Conventions}\label{Sec.Results.Notations}
 We denote by $\ux$ the multivariable $x_1,\dots,x_p$.

Let $A\in\Mat$ and denote by $I_j(A)$ the determinantal ideal of all the $j\times j$ minors.
Denote by $ann.coker(A)$ the
annihilator-of-cokernel ideal of the homomorphism of free modules
$R^{\oplus n}\stackrel{A}{\to}R^{\oplus m}$. The
properties/relation/computability of these ideals are discussed in
\S\ref{Sec.Background.Fitting.Ideals.ann.coker}.

We often use the ideal quotient, $(I:J)=\{f\in R|\ fJ\sseteq I\}$ and the integral closure of ideals $I\sseteq\overline{I}$,
 see \S\ref{Sec.Background.Integral.Closure.Ideals}.
  Sometimes we write $I:J$ instead of $(I:J)$, to avoid multiple-brackets formulas.

Sometimes we take the $\cm$-adic completion, $R\to\hR$, this induces  the completion map $\Mat\to Mat_{m\times n}(\hR)$,
denote by $\hA$ the image of $A$.

\

Suppose $R$ is local and $J\supseteq\cm^\infty$. The Loewy length, $ll_R(J)\le\infty$,  is the minimal number $N\le\infty$
 that satisfies: $J\supseteq\cm^N$. This number also equals the degree of the socle of the quotient module $\quots{R}{J}$. For $R=\k[[\ux]]$ we have yet another expression, via the Castelnuovo-Mumford regularity,
 $ll_R(J)=reg(\quots{R}{J})+1$, see \cite[exercise 20.18]{Eisenbud}.

\

For the applications to the (in)finite determinacy we assume (in such cases we state this explicitly):
\beq\label{Eq.Intro.assumption.R.Noetherian.or.surjects}
\ber
\text{$(R,\cm)$ is a local ring over a field $\k$ of zero characteristic, and }
\\
\text{either $R$ is Noetherian or  the $\cm$-adic completion map is surjective, $R\twoheadrightarrow\hR$, and $\hR$ is Noetherian.}
\eer\eeq
In this case by $dim(R)$ we mean $dim(\hR)$.

Many rings satisfy this condition, e.g. $\quots{\k[\ux]}{I}$, $\quots{\k\langle \ux\rangle}{I}$ (algebraic power series), $\quots{\k\{\ux\}}{I}$ (convergent power series),
 $\quots{\k[[\ux]]}{I}$ and $C^\infty(\R^p,0)$. (For the later see  see \S\ref{Sec.Background.Approximation.Properties}.)
But the ring $C^r(\R^p,0)$, of r-times continuously differentiable germs, for $0<r<\infty$,
 does not satisfy condition \eqref{Eq.Intro.assumption.R.Noetherian.or.surjects}.  Its completion is not Noetherian, e.g. because the ideal
$\cm^{r+0^+}:=\{g\cdot\cm^r|\ g\in C^0(\R^p,0),\ g(0)=0\}$
is not finitely generated, and neither is its image in $\widehat{C^r(\R^p,0)}$.

\

In some cases we assume that the ring $R$ is a principal ideal domain (PID) or a discrete valuation ring (DVR) or a unique factorization domain (UFD),
 then we state this explicitly.

\subsubsection{Generic properties}\label{Sec.Intro.Genericity.Tougeron}
In this paper the {\em genericity} is always used in the following sense of \cite{Tougeron1968}. Fix some ``parameter space" $M$,
 typically a finitely generated $R$-module, suppose the $jet_d$-projections $\{R\stackrel{jet_d}{\to}\quots{R}{\cm^{d+1}}\}_{d\in\N}$ induce
 the projections $M\stackrel{jet_d}{\to}jet_d(M)$. We assume that all $jet_d(M)$ are  algebraic subschemes over $\k$.
A property $\cP$ is said to hold generically in $M$ if it holds for all points of $M$ lying in the complement to a
 subset $X\sset M$ which is ``of infinite codimension".
Namely, for any $d$ the projection of $X$ lies inside some algebraic subschemes, $jet_d(X)\sseteq Y_d\sset jet_d(M)$, such that $codim_{jet_d(M)}(Y_d)\to\infty$.

\

\bed\label{Def.Generic.Finite.Determinacy}
We say that the generic finite determinacy holds for a given action $G\circlearrowright\Si$ if for any $A\in\Si$, any number $N<\infty$ and
 the generic matrix $B\in Mat_{m\times n}(\cm^N)$, such that $A+B\in\Si$, the matrix $A+B$ is finitely determined.
\eed

Equivalently: the $\cm$-adic closure of the set of finitely-determined matrices is the whole space $\Si$.

Recall that in singularity theory the dimension of the vector space $T^1_{(\Si,G,z)}$ is called the Tjurina number, for a given deformation scenario $(\Si,G,A)$.
 The finite determinacy means the finiteness of the relevant Tjurina number.
 The genericity of finite determinacy  means:
 the stratum of matrices whose Tjurina number is infinite, $\Si^G_{\tau=\infty}\sset \Si$, is itself of infinite codimension.

\subsection{The criteria for the $G_r$, $G_l$, $G_{lr}$ actions of Example \ref{Ex.Intro.Typical.Groups}}\label{Sec.Results.Criteria.Glr}
\bthe\label{Thm.Results.Annihilators.for.GlGrGlr} Let $A\in\Si=\Mat$, with $m\le n$.

\bee
\item $T^1_{(\Si,G_r,A)}=\Big(coker(A)\Big)^{\oplus n}$ and therefore $ann(T^1_{(\Si,G_r,A)})=ann.coker(A)$.
\item If $m<n$ then $ann\big(T^1_{(\Si,G_l,A)}\big)=\{0\}$. If $m=n$ then $ann\big(T^1_{(\Si,G_l,A)}\big)=ann.coker(A)$.
\item $ann\big(T^1_{(\Si,G_{lr},A)}\big)\supseteq ann.coker(A)$.
 If $R$ is Noetherian and $I_{m-1}(A)\neq\{0\}$ then
 \[\overline{ann\big(T^1_{(\Si,G_{lr},A)}\big)}\sseteq \overline{I_m(A)}:\overline{I_{m-1}(A)}.
 \]
\eee
\ethe
The bounds in part 3 are close and in many cases  we have
 $ann\big(T^1_{(\Si,G_{lr},A)}\big)=ann.coker(A)$, see examples in \S\ref{Sec.Proofs.Glr}.

 In fact, in \S\ref{Sec.Proofs.Glr} we prove a slightly stronger statement:
\beq
\overline{ann\big(T^1_{(\Si,G_{lr},A)}\big)\cdot  I_{m-1}(A)}\sseteq \overline{ I_m(A)}.
\eeq


\

Theorem \ref{Thm.Results.Annihilators.for.GlGrGlr} bounds $ann\big(T^1_{(\Si,G_{lr},A)}\big)$ in terms of
$ann.coker(A)$ and $\overline{I_m(A)}:\overline{I_{m-1}(A)}$. These ideals are usually rather small,
thus finite-$G_{lr}$-determinacy places severe restrictions on the ring $R$.

\bprop\label{Thm.Results.Fin.Det.for.GlGrGlr}  Suppose $R$ satisfies condition \eqref{Eq.Intro.assumption.R.Noetherian.or.surjects} and let $\Si=\Mat$.
\bee
\item
If $ann.coker(A)\supseteq\cm^\infty$ then $ll_R\Big(ann.coker(A)\Big)-1 \le ord^\Si_{G_r}(A)\le ll_R\Big(ann.coker(A|_{\cm R^{\oplus n}})\Big)-1.$

 In particular:
\bee[i.]
\item If $dim(R)>(n-m+1)$ then no matrix in $\Matm$ is finitely-$G_r$-determined.
\item If $dim(R)\le(n-m+1)$ then the generic $G_r$-finite determinacy holds.
\eee
\item If $m<n$ then no matrix $A\in\Matm$ is finitely-$G_l$-determined.
 If $m=n$ then  $ord^\Si_{G_r}(A)=ord^\Si_{G_l}(A)$.
\item ($G_{lr}$-equivalence)
\bee[i.]
\item
If $dim(R)>(n-m+1)$ then no matrix in $\Matm$ is finitely-$G_{lr}$-determined,
\item  Suppose $ann.coker(A)\supseteq\cm^\infty$ then
\[ll_\hR\Big(\overline{I_m(\hA)}:\overline{I_{m-1}(\hA)}\Big)-1\le ord^\Si_{G_{lr}}(A)\le ll_R\Big(ann.coker(A|_{\cm\cdot R^{\oplus n}})\Big)-1.\]
\eee
\eee
\eprop
Here  $ann.coker(A|_{\cm R^{\oplus n}})$ denotes the annihilator-of-cokernel of the restricted map $\cm\cdot R^{\oplus n}\stackrel{A}{\to}R^{\oplus m}$.

Regarding the notations in part 3.ii:
 $\hA$ denotes the image of $A$ under the completion map $R\to\hR$.  If $R$ is Noetherian then
 \beq
 ll_\hR\Big(\overline{I_m(\hA)}:\overline{I_{m-1}(\hA)}\Big)=ll_R\Big(\overline{I_m(A)}:\overline{I_{m-1}(A)}\Big),
 \eeq
  see Corollary \ref{Thm.Loewy.Length.under.completion}, so this number can be computed in either way.

\

The bounds in this proposition are rather tight, e.g.:
\beq
ll_R\big(ann.coker(A)\big)\le ll_R\big(ann.coker(A|_{\cm R^{\oplus n}})\big)\le ll_R\big(ann.coker(A)\big)+1
\eeq
 and
 also the ideals $\overline{I_m(A)}:\overline{I_{m-1}(A)}$, $ann.coker(A)$ are quite close and often coincide.
 The proofs, examples and other corollaries are given in \S\ref{Sec.Proofs.Glr}.

\

We remark that part 1 of the proposition implies also criteria for the infinite-$G_r$-determinacy, see Corollary \ref{Thm.Proofs.Infinite.Determinacy.Gr} for the $C^\infty(\R^p,0)$ case.

\subsection{The criteria for $G_{congr}$}\label{Sec.Results.Congruence}
Through this subsection we assume $m>1$.

Recall that for skew-symmetric matrices of even size holds $I_m(A)=(det(A))=(Pf(A))^2$, the square of the Pfaffian ideal of $A$. Moreover, $ann.coker(A)\supseteq Pf(A)$, see Proposition \ref{Thm.Background.AntiSymmetric.Matrices.Pfaffians}.

For skew-symmetric matrices of odd size: $det(A)=0$ and $ann.coker(A)=\{0\}$. As the measure of non-degeneracy of such matrices we use the $(m-1)$'st Pfaffian ideal
\beq
(Pf_{m-1}(A)):=\suml_{i=1}^m (Pf(A_{\hi})).
\eeq
Here $A_{\hi}$ is the $(m-1)\times(m-1)$ block of $A$, obtained by erasing the $i$'th row and column.
(Note that $A_{\hi}$ is skew-symmetric and $(Pf(A_{\hi}))^2=I_{m-1}(A_{\hi})$.)
See \S\ref{Sec.Background.Pfaffians} for more detail.
\bthe\label{Thm.Results.Annihilator.T1.for.Gcongr}  Suppose $I_{m-1}(A)\neq\{0\}$.

\bee
\item Let $\Si=Mat_{m\times m}(R)$, suppose $R$ satisfies the condition \eqref{Eq.Intro.assumption.R.Noetherian.or.surjects} and
 $dim(R)>0$.
 Then  $ann(T^1_{(\Si,G_{congr},A)})\sseteq \cm^\infty$.
\\ In particular, if $\cm^\infty=\{0\}$ then $ann(T^1_{(\Si,G_{congr},A)})=\{0\}$.

\item Let $A\in \Si=Mat^{sym}_{m\times m}(R)$ then $ann.coker(A)\sseteq ann(T^1_{\Si,G_{congr},A)})$.
\\ Suppose $R$ is Noetherian then $ann(T^1_{\Si,G_{congr},A)})\sseteq\overline{I_{m}(A)}:\overline{I_{m-1}(A)}$.

\item
\bee[i.]
\item Let $A\in \Si=Mat^{skew-sym}_{m\times m}(R)$ with $m$ even, then
$ann.coker(A)\sseteq ann(T^1_{(\Si,G_{congr},A)})$.
\\ If  $R$ is Noetherian  then $ann(T^1_{(\Si,G_{congr},A)})\sseteq\overline{I_{m}(A)}:\overline{I_{m-1}(A)}$.

\item  Let $A\in \Si=Mat^{skew-sym}_{m\times m}(R)$ with $m$ odd, then
$(Pf_{m-1}(A))\sseteq ann(T^1_{(\Si,G_{congr},A)})$.
\\ If  $R$ is Noetherian  then $ann(T^1_{(\Si,G_{congr},A)})\sseteq\overline{I_{m-1}(A)}:\overline{I_{m-2}(A)}$.
\eee
\eee
\ethe
The bounds in parts 2,3 are rather close and often coincide, see examples in \S\ref{Sec.Proofs.Congruence}.

\

The bounds on $ann(T^1_{(\Si,G_{congr},A)})$ give immediate applications to the finite determinacy:

\bprop\label{Thm.Results.Finite.Determinacy.G_congr}
Suppose $R$ satisfies the condition \eqref{Eq.Intro.assumption.R.Noetherian.or.surjects}.

{\it (1)} If $dim(R)>1$ then no $A\in Mat^{sym}_{m\times m}(\cm)$ is finitely-$(Mat^{sym}_{m\times m}(R),G_{congr})$-determined.

Similarly, no $A\in Mat^{skew-sym}_{m\times m}(\cm)$ with $m$-even is finitely-$(Mat^{skew-sym}_{m\times m}(R),G_{congr})$-determined.

{\it (1')} If $m$ is odd and $dim(R)>3$ then no $A\in Mat^{skew-sym}_{m\times m}(\cm)$ is finitely-$(Mat^{skew-sym}_{m\times m}(R),G_{congr})$-determined.

{\it (1'')} If $dim(R)>0$ then no $A\in Mat_{m\times m}(\cm)$ is finitely-$(Mat_{m\times m}(R),G_{congr})$-determined.

{\it (2)} Suppose $R$ is Henselian.

\quad {\bf i.} Suppose $dim(R)=1$ and either $\Si=Mat^{sym}_{m\times m}(R)$ or $\Si=Mat^{skew-sym}_{m\times m}(R)$ for $m$-even.
Let $A\in\Si$ and suppose $I_m(A)\supsetneq\cm^\infty$. Then
\[
ll_\hR\Big(\overline{I_m(\hA)}:\overline{I_{m-1}(\hA)}\Big)-1\le ord^{\Si}_{G_{congr}}(A)\le ll_R(ann.coker(A)).
\]
 In particular, $A$ is finitely-$(\Si,G_{congr})$-determined iff $det(A)$ is not a zero divisor in $R$.
 In particular, the generic finite determinacy holds.

\quad {\bf ii.} Suppose $dim(R)\le3$ and $A\in \Si=Mat^{skew-sym}_{m\times m}(R)$, for $m$-odd.
Suppose moreover $I_{m-1}(A)\supsetneq\cm^\infty$.
 Then
\[ll_\hR\Big(\overline{I_{m-1}(\hA)}:\overline{I_{m-2}(\hA)}\Big)-1\le ord^{\Si}_{G_{congr}}(A)\le ll_R(Pf_{m-1}(A)).
\]
In particular, $A$ is finitely-$(\Si,G_{congr})$-determined iff either $I_{m-1}(A)=R$ or
$I_{m-1}(A)\sset R$ is of height$=dim(R)$, as expected.
 In particular, the generic finite determinacy holds.
\eprop

The proofs, examples and further corollaries are in \S\ref{Sec.Proofs.Congruence}.

 For the criteria of infinite
determinacy of $C^\infty(\R^p,0)$-valued (skew-)symmetric forms see Corollary \ref{Thm.Proofs.Infinite.Determinacy.G_congr}.

\subsection{Upper-block-triangular matrices.}
For  an upper-block-triangular matrix $A\in Mat^{up}_{\{m_i\}\times\{n_j\}}(R)$ take the block-structure, $\{A_{ij}\}_{1\le i,j\le k}$,
 as in Example \ref{Ex.Intro.Typical.Groups}.
The context/motivation  for this case, the definitions of $T^1_{(\Si,G^{up}_r,A)}$, $T^1_{(\Si,G^{up}_{lr},A)}$,
their roles, the proofs, examples  and corollaries  are in \S\ref{Sec.Proofs.Upper.Triangular}.

 Denote by   $\cornA{q}$ the sub-matrix of $A$ consisting of the
 first $\suml^q_{i=1}m_i$ rows and the first $\suml^q_{j=1}n_j$ columns. (Thus $\cornA{1}$ is the first upper-left block, while $\cornA{k}=A$.)
 We use this notation in part (1) of the theorem:
\bthe\label{Thm.Results.Annihilator.T1.Upper.Triang}
Let $A\in \Si=Mat^{up}_{\{m_i\}\times\{n_j\}}(R)$.
\bee
\item $ann(T^1_{(\Si,G^{up}_r,A)})=\capl^k_{q=1}ann.coker(\cornA{q})$.
\item
In particular,
 $\prodl_{q=1}^k ann.coker(A_{qq})\sseteq ann(T^1_{(\Si,G^{up}_r,A)})\sseteq\capl_{q=1}^k ann.coker(A_{qq})$.
\item Suppose $R$ is Noetherian then
\[
\prodl_{q=1}^k ann.coker(A_{qq})\sseteq \capl^k_{q=1}ann.coker(\cornA{q})\sseteq ann(T^1_{(\Si,G^{up}_{lr},A)})\sseteq\capl_{q=1}^k
 \Big(\overline{I_{m_q}(A_{qq})}:\overline{I_{m_q-1}(A_{qq})}\Big).
\]
\eee
\ethe
In part 2 both bounds are sharp, see examples in \S\ref{Sec.Proofs.Upper.Triangular}.

\subsection{The non-finite determinacy for conjugation}\label{Sec.Results.Conjugation}
 The conjugation action is defined in Example \ref{Ex.Intro.Typical.Groups}.
\bthe\label{Thm.Results.T1.Conjugation} Suppose $R$ satisfies condition \eqref{Eq.Intro.assumption.R.Noetherian.or.surjects} and let $A\in Mat_{m\times m}(R)$.
\bee
\item $ann(T^1_{(Mat_{m\times m}(R),G_{conj},A)})\sseteq\cm^\infty$.
\item In particular, if $dim(R)>0$ then there are no finitely-$G_{conj}$-determined matrices.
\eee
\ethe

\subsection{Finite determinacy of chains of free modules}
As an immediate application of our methods we consider bounded chains (or complexes) of free modules.
Take such a chain, $\cdots\stackrel{\phi_{i+1}}{\to} F_i\stackrel{\phi_i}{\to}F_{i-1}\stackrel{\phi_{i-1}}{\to}\cdots$,
fix some  bases of $\{F_i\}$, so that the maps are represented by some matrices. The deformations of such a chain are taken up to isomorphisms, induced by the action of the product of groups $\prodl_i GL(F_i,R)$. We obtain the bounds on the determinacy in  \S\ref{Sec.Proofs.Complexes}, quantifying and generalizing theorem 5.8 of \cite{Cutkosky-Srinivasan}.

\subsection{Algebraization }\label{Sec.Results.Algebraization}
 Suppose $\quots{\k[\ux]}{I}\sseteq R$   and the completion of $R$ at the ideal $(\ux)$ is $\hR=\quots{\k[[\ux]]}{I}$.
  For example, this holds when  $R$ is a local ring over a field $\k$, by Cohen Structure Theorem.

A trivial consequence of finite determinacy is the  algebraization:
the  $G$-orbit of  a given object contains  algebraic representatives (with entries in $\quots{\k[\ux]}{I}$)
 and they are obtained by just ``cutting the
 Taylor-tails". Moreover, the order of determinacy gives an upper bound on the possible degrees of polynomials.
 As we show in \S\ref{Sec.Background.Fin.Determin.vs.Stable.Algebraizability}, finite determinacy is equivalent to {\em stable algebraizability:}
   any higher order deformation of
  a given object is $G$-equivalent to an algebraic family.

Even if $A$ is not finitely determined, i.e. $ann(T^1_{(\Si,G,A)})$ does not contain any $\cm^N$, one can consider
some ideal $J$ such that $J+ann(T^1_{(\Si,G,A)})$ contains a power of the maximal ideal. Then one gets the ``algebraization modulo $J$", or the ``algebraization with respect to some subset of variables". Therefore our results generalize those of \cite{Elkik73}, \cite{Kucharz} and many others.

\subsection{Relative determinacy/admissible deformations}\label{Sec.Results.Relative.Determinacy}
As we see in Theorems \ref{Thm.Results.Annihilators.for.GlGrGlr}, \ref{Thm.Results.Annihilator.T1.for.Gcongr}, \ref{Thm.Results.Annihilator.T1.Upper.Triang}, the condition $ann(T^1_{(\Si,G,A)})\supseteq\cm^N$ can be quite restrictive. For example:
\bei
\item for the action $G_{lr}\circlearrowright\Mat$ it means:
$ann.coker(A)\supseteq\cm^N$, i.e. the module $coker(A)$ is of finite length. Geometrically: the module is supported only at one point, the origin, i.e. is a skyscraper.
\item
for the action $G_{congr}\circlearrowright Mat^{sym}_{m\times m}(R)$ it means:  $A$ and its associated quadratic form are non-degenerate
 on the punctured neighborhood of the origin.
\eei
  However, one often needs to deform $A$ not inside the whole $\Mat$, or $\Si$, but only by elements of $Mat_{m\times n}(J)$ (for some ideal $J\sset R$),
   i.e. the deformation space is
  $\Si^{(J)}:=\Si\cap\Big(\{A\}+Mat_{m\times n}(J)\Big)$.
 For example, one studies the determinacy for the action $G_{lr}\circlearrowright Mat_{m\times n}(J)$ or $G_{congr}\circlearrowright Mat^{sym}_{m\times m}(J)$.

  Similarly, starting from a group action, $G\circlearrowright\Mat$, one restricts to the transformations preserving the filtration on $\Mat$.
   Often one restricts further, to the subgroup of transformations that are identities modulo some (prescribed) ideal $I\subsetneq R$. More precisely:
\beq\label{Def.Unipotent.Subgroup}
G^{(I)}:=\Big\{g\in G|\ g\cdot Mat_{m\times n}(I)=Mat_{m\times n}(I)\ and\ [g]=[Id]\circlearrowright\quotient{\Mat}{Mat_{m\times n}(I)}\Big\}.
\eeq
For example, for $I=\cm^{k}$ the group $G^{(\cm^k)}$ consists of elements that are identities up to the order
 $(k-1)$.

\

Working with $\Si^{(J)}$ and $G^{(I)}$  leads to the notion of relative determinacy/admissible deformations.
For complex analytic non-isolated hypersurface singularities this was studied in \cite{Pellikaan88}, \cite{Pellikaan}, \cite{Siersma83}, \cite{Siersma},
\cite{de Jong-van Straten1990}, \cite{Jiang}. 
 For the $C^\infty(\R^p,0)$-version see \cite{Grandjean00}, \cite{Grandjean04}, \cite{Thilliez}, \cite{Cutkosky-Srinivasan}.

Our criteria directly adapt to this formulation, see Proposition \ref{Thm.Proofs.Corollary.For.Rel.Determinacy} in \S\ref{Sec.Proofs.Relative.Determinacy}.

\subsection{Remarks}\label{Sec.Results.Remarks}
\subsubsection{}
 Theorems  \ref{Thm.Results.Annihilators.for.GlGrGlr}, \ref{Thm.Results.Annihilator.T1.for.Gcongr}, \ref{Thm.Results.Annihilator.T1.Upper.Triang}
 and their corollaries
 are absolutely explicit, ready-to-use criteria. Their role is similar to the classical Mather theorem  (for $\k=\bar\k$), the transition from:

$\bullet$ ``the hypersurface germ is finitely determined iff its miniversal deformation is finite dimensional"
\\to

 $\bullet$ ``the hypersurface germ is finitely determined iff it has at most an isolated singularity".

Moreover, we do not just establish the finite determinacy, but give rather tight bounds for the order of determinacy.

To emphasize, in most cases the module  $T^1_{(\Si,G,A)}$ does not admit any explicit/short presentation, e.g. for  $T^1_{(\Mat,G_{lr},A)}$
 one cannot write anything more explicit than $Ext^1_R(coker(A),coker(A))$. But the bounds we get are quite
 explicit and  tight, the upper and lower bounds often have the same integral closure.

\subsubsection{``Restrictiveness" of the results} As our results read, for many actions there are no finitely determined matrices in $\Matm$.
For example, there are no finitely determined $\cm$-valued quadratic
forms when $dim(R)>1$ (Proposition
\ref{Thm.Results.Finite.Determinacy.G_congr}). For the action
$G_{conj}\circlearrowright Mat_{m\times m}(R)$ the situation is even
worse, Theorem \ref{Thm.Results.T1.Conjugation}.

Sometimes one insists on having enough finitely determined objects, e.g. one wants the generic finite determinacy in the sense of definition \ref{Def.Generic.Finite.Determinacy}. The two ways to achieve this are:
\bei
\item Either to enlarge the group, e.g. to use also the local coordinate changes/automorphisms
 of the ring, $Aut_\k(R)$. Then one considers the action of $G_{lr}\rtimes Aut_\k(R)$, $G_{congr}\rtimes Aut_\k(R)$, etc. We do this in \cite{Belitski-Kerner2}.
\item Or to restrict the possible deformations, e.g. to deform in a way that preserves the determinantal ideals, $\{I_j(A)\}$, or to preserve
 the characteristic polynomial for square matrices. We do this in \cite{Belitski-Kerner3}.
\eei

\subsubsection{A ``theoretical" remark}
 Fix some action $G\circlearrowright\Si\sseteq Mat_{m\times n}(R)$, we assume that $T_{(\Si,A)}$ is a free $R$-module. (This holds in our cases.)
  Take the generating matrix $\cA_{G,A}$ of the module $T_{(GA,A)}$ so that $Im(\cA_{G,A})=T_{(GA,A)}\sseteq T_{(\Si,A)}$. Suppose $\cA_{G,A}$ is of size $m_1\times n_1$, here $m_1=rank(T_{(\Si,A)})$, then
\beq
T^1_{(\Si,G,A)}=\quotient{T_{(\Si,A)}}{T_{(GA,A)}}=
coker(\cA_{G,A}),\quad
\text{while}\quad T^1_{(Mat_{m_1\times n_1}(R),G_r,\cA_{G,A})}=coker(\cA_{G,A})^{\oplus n_1}.
\eeq
(Here on the right $G_r=GL(n_1,R)$.) In particular, $ann(T^1_{(\Si,G,A)})=ann(T^1_{(Mat_{m_1\times n_1}(R),G_r,\cA_{G,A})})$.

Therefore the study of $T^1$ for any $G$ action  is ``embedded" into the study of
 $T^1$ for the $G_r$ action, on the space of bigger size matrices.

\subsubsection{Finite determinacy of (non-)embedded modules and (skew-)symmetric forms}
The notion of finite determinacy and our results can be stated more algebraically as follows.

 The projection $R\stackrel{jet_k}{\to}\quots{R}{\cm^{k+1}}$ induces the functor, $mod(R)\to mod({\quots{R}{\cm^{k+1}}})$,
defined by $M\to \quots{M}{\cm^{k+1} M}$. This functor is trivially
surjective on the objects, by the restriction of scalars.  But the
functor is not injective on the objects, e.g. if the presentation
matrix of $M$ has all its entries in $\cm^{k+1}$ then
$\quots{M}{\cm^{k+1} M}$ is a free $\quots{R}{\cm^{k+1}}$ module.
Our results give the regions of parameters (the size of presentation
matrix, $dim(R)$) for which the functor is ``generically injective on
the objects".

\

 Let $M\sseteq R^{\oplus m}$ be a finitely
generated $R$-submodule. Fix some set of generators of $M$, combine them into the matrix $A$.
\bei
\item
The  finite-$G_r$-determinacy of $A$ means:
 ``{\em if $M,N\sseteq R^{\oplus m}$ and $jet_k M=jet_k N$ for $R\stackrel{jet_k}{\to}\quots{R}{\cm^{k+1}}$, then $M=N\sseteq R^{\oplus m}$}".
\item The finite-$G_{lr}$-determinacy of $A$ means:
 ``{\em if $M,N\sseteq R^{\oplus m}$ and $jet_k M=jet_k N$ for $R\stackrel{jet_k}{\to}\quots{R}{\cm^{k+1}}$, then $M= \phi(N)\sseteq R^{\oplus m}$ for some $\phi\in GL(m,R)$}".
\eei

\

Similarly for (skew-)symmetric forms the finite determinacy means that the form is determined (up to the change of generators) by its $k$'th jet for $k\gg0$.

\subsection{Relation to Singularity Theory}\label{Sec.Results.Relation.to.Singularity.Theory}

The (in)finite determinacy is a classical notion of Singularity Theory, well studied for functions and maps over the ``classical" rings,
 when $\k=\R$ or $\C$, and $R=\k[[\ux]]$, $\k\{\ux\}$, $\k[\ux]_{(\ux)}$, $C^{\infty}(\R^p,0)$,
\cite{Mather1968}, \cite{Wilson81}, \cite{Wilson82}, \cite{Damon1984}, \cite{Wall-1981}, \cite{AGLV}, etc.

The generic finite determinacy has been considered in \cite{Wall-1979}, \cite{du Plessis1982}, \cite{du Plessis1983}.

(For a short discussion of the results on the determinacy of functions/maps see \cite[\S2]{Belitski-Kerner}.)

\subsubsection{Some known results on ``matrix singularities"}\label{Sec.Results.Previous.Results.Singular.Theory}
Some necessary and sufficient conditions for  finite determinacy of modules over the ring $C^\infty(\R^p,0)$ were obtained in \cite{Tougeron1968}.

The square matrices over $R=\k\{x_1,\dots,x_p\}$, for $\k\in\R,\C$,  and $G=G_{lr}\rtimes Aut(R)$,
were considered in \cite{Haslinger}, \cite{Bruce2003}, \cite{Bruce-Tari04}, and further studied in \cite{Bruce-Goryun-Zakal02},
\cite{Bruce2003}, \cite{Goryun-Mond05}, \cite{Goryun-Zakal03}. In particular, the generic finite
determinacy was established and the simple types were classified.
 We mention also the recent works,  \cite{Pereira.phd}, \cite{Pereira.2017}, \cite{Ahmed-Ruas}.

In \cite{Cutkosky-Srinivasan} they study the finite determinacy of matrices for the actions of $G_r$, $G_{lr}$, $G_r\rtimes Aut_\k(R)$, $G_{lr}\rtimes Aut_\k(R)$ on matrices over complete local rings.
 They obtain the qualitative results for the finite determinacy.

In \cite{Damon.Pike.1}, \cite{Damon.Pike.2} they study the vanishing topology of matrix singularities, relating these to the free divisors. In \cite{Fr.Kr-Za.2015} they study the vanishing topology of determinantal singularities.

The study of matrix determinacy in positive characteristic has been initiated in \cite{Greuel-Pham}, \cite{Pham}.

\

We emphasize that most of the previous works addressed only the general criteria for the finite determinacy. In our work we directly compute (or at least bound) the order(s) of determinacy.

\subsubsection{}
In our results we see the standard dichotomy: either the finite determinacy is
the generic property  or there are no finitely determined matrices in $\Matm\cap\Si$.
This is the analogue of the ``nice/bad dimensions" for the determinacy of maps $Maps\Big((\k^n,0),(\k^m,0)\Big)$,
\cite[III.1]{AGLV}, \cite{Mather1968}.

\section{Background and Preparations}\label{Sec.Background}
 In this section we collect the standard facts of commutative algebra, thus, unless stated otherwise,  $R$ is just a commutative unital ring.
We denote the zero matrix by $\zero$, the identity matrix by $\one$.

\
Let $(R,\cm)$ be a local ring, with residue field $\k=\quots{R}{\cm}$.
Given $A\in Mat_{m\times n}(R)$ we sometimes  ``compute it at the origin", i.e. take its image over the residue field, $A|_0\in Mat_{m\times n}(\k)$.
 In particular, $A|_0=\zero$ means $A\in Mat_{m\times n}(\cm)$.  As $A|_0$ is a matrix over a field, we take its classical rank, $r=rank(A|_0)$.
Recall that the rank of a skew-symmetric matrix is necessarily even.

\

 We denote by $E=\oplusl_i \bpm 0&1\\-1&0\epm$ the canonical skew-symmetric matrix,
  its size (i.e. the number of summands) should be clear from the  context.

\subsection{The relevant canonical forms for matrices}\label{Sec.Background.Canon.Forms.Matrices}

\bprop\label{Thm.Background.Chip.off.Unit} Let $(R,\cm)$ be a local ring.
\bee
\item Let $A\in\Mat$ then $A\stackrel{G_{lr}}{\sim}\one_{r\times r}\oplus\tA$, where $r=rank(A|_0)$ and
 $\tA\in Mat_{(m-r)\times(n-r)}(\cm)$.
\item Let $A\in Mat^{sym}_{m\times m}(R)$ then $A\stackrel{G_{congr}}{\sim}(\oplusl_i \la_i\one_i)\oplus\tA$,
 where $\{ \la_i\in R\smin\cm\}_i$ and $\tA\in Mat^{sym}_{(m-r)\times(m-r)}(\cm)$.
\item Let $A\in Mat^{skew-sym}_{m\times m}(R)$ then $A\stackrel{G_{congr}}{\sim}\Big(\oplusl_{i} \la_i E_i\Big)\oplus\tA$,
 where $\{ \la_i\in R\smin\cm\}_i$  and $\tA\in Mat^{skew-sym}_{(m-r)\times(m-r)}(\cm)$.
\eee
\eprop
Part 1 is known in Commutative Algebra e.g. as ``passing to the minimal resolution of the module $coker(A)$".
 Similarly, in parts 2,3 we split a (skew-)symmetric form into its regular part (in the canonical form) and the complement - the `purely degenerate' part.

(Part 2. is proved e.g. in \cite[theorem 3, page 345]{Birkhoff-MacLane}, for part 3 see \cite[exercise 9, page 347]{Birkhoff-MacLane}.)

\

Let $R$ be a discrete valuation ring (DVR) over $\k$. In particular $R$ is local, regular and $dim(R)=1$.
The simplest examples of such rings are the formal power series, $R=\k[[t]]$, or the analytic series, $R=\k\{t\}$, when $\k$ is a normed field.
 Note that the ring of real-valued smooth function-germs, $C^\infty(\R^1,0)$, is not a DVR.
Matrices over DVR often have good canonical forms.

\bprop\label{Thm.Background.Canonical.Form.Matrix.over.DVR} Suppose $R$ is a DVR over $\k$.
\bee
\item For any $A\in \Mat$ holds:
 $A\stackrel{G_{lr}}{\sim}\bpm \la_1&0&\dots&\dots&0\\0&\la_2&0&\dots&\dots\\\dots&\dots&\dots&\dots&\dots\\0&\dots&&\la_{m}&0\epm$,
 where $R\supseteq(\la_1)\supseteq(\la_2)\supseteq\cdots\supseteq(\la_m)$.
\item For any $A\in  Mat^{sym}_{m\times m}(\cm)$ holds: $A\stackrel{G_{congr}}{\sim}\oplusl_i \la_i\one_i$,
 where $R\supseteq(\la_1)\supseteq(\la_2)\supseteq\cdots\supseteq(\la_m)$.
\item For any $A\in Mat^{skew-sym}_{m\times m}(\cm)$ holds: $A\stackrel{G_{congr}}{\sim}\Big(\oplusl_{i} \la_i E_i\Big)\oplus \zero$, where $R\supseteq(\la_1)\supseteq(\la_2)\supseteq\cdots$, while $\zero$ is the zero matrix of an appropriate size.
\eee
\eprop
Here part 1 is the Smith normal form.
Considering $A$ as a presentation matrix of the $R$-module $coker(A)$, the statement is:
 ``every module over a DVR is a direct sum of cyclic modules". Parts 2,3 read: every (skew-)symmetric form over a DVR splits into the
  direct sum of rank-one forms.

While parts 2,3 are standard facts, we did not find an exact reference, thus we sketch a proof. Let $t$ be a uniformizing parameter
of $R$. Denote by $p$ the minimal among the vanishing orders of the entries of $A$. Then $t^{-p}A$ is a (skew-)symmetric matrix over
$R$ and $t^{-p}A|_0\neq\zero$. Apply Proposition \ref{Thm.Background.Chip.off.Unit} to get,  in the symmetric case,
 $A\stackrel{G_{congr}}{\sim}t^{p}\Big((\oplus_i\la_i\one_i) \oplus \tA\Big)$,
 (or $A\stackrel{G_{congr}}{\sim}t^{p}\Big((\oplus_i\la_i E_i) \oplus \tA\Big)$ in the skew-symmetric case).
 Thus we get the reduction in size and one repeats for $\tA$.

\


\subsection{Determinantal ideals and the annihilator-of-cokernel}\label{Sec.Background.Fitting.Ideals.ann.coker}
\cite[\S 20]{Eisenbud}
For $0< j\le m$ and $A\in\Mat$ denote by $I_j(A)\sset R$
the $j$'th determinantal ideal, generated by all the $j\times j$ minors of $A$.
By definition $I_{j\le 0}(A)=R$ and $I_{j>m}(A)=\{0\}$. When the size of $A$ is not stated explicitly we denote the ideal of maximal minors by $I_{max}(A)$.
 The chain of ideals
 $R=I_0(A)\supseteq I_1(A)\supseteq\cdots\supseteq I_m(A)$ is {\em invariant} under the $G_{lr}$-action.

Determinantal ideals behave nicely under ring homomorphisms,  for any $R\stackrel{\phi}{\to}S$ and the induced map
 $\Mat\stackrel{\phi}{\to}Mat_{m\times n}(S)$ one has: $I_j(\phi(A))=\phi(I_j(A))$.

\subsubsection{The expected heights} (These will be used later for the generic finite determinacy/bounds on the annihilator.)

If $A|_0=\zero$, i.e. $A\in\Matm$, then
 the height of $I_j(A)$ is at most $min\Big((m+1-j)(n+1-j),dim(R)\Big)$ and this bound is achieved generically.
 More precisely,  for any   $A\in\Matm$, any $N>0$ and any generic $B\in Mat_{m\times n}(\cm^N)$ the height of the ideal $I_j(A+B)$ equals
 $min\Big((m+1-j)(n+1-j),dim(R)\Big)$.

For $A\in Mat^{sym}_{m\times m}(\cm)$ the expected height of $I_j(A)$ is $min\Big(\bin{m-j+2}{2},dim(R)\Big)$. For $A\in Mat^{skew-sym}_{m\times m}(\cm)$  and $j$  even the expected height of $I_j(A)$ is
$min\Big(\bin{m-j+2}{2},dim(R)\Big)$.
 For $A\in Mat^{skew-sym}_{m\times m}(\cm)$  and $j$  odd the expected height of $I_j(A)$ is
 $min\Big(\bin{m-j+1}{2},dim(R)\Big)$.

\subsubsection{$ann.coker(A)$}

A matrix $A\in\Mat$ can be considered as a homomorphism of free modules, its cokernel is an $R$-module as well:
\beq
 R^{\oplus n}\stackrel{A}{\to} R^{\oplus m}\to coker(A)\to0.
 \eeq
  Then one takes the annihilator-of-cokernel ideal:
\beq
ann.coker(A)=\{f\in R|\ f\cdot R^{\oplus m}
\sseteq
Im(A)\}=ann\quots{R^{\oplus m}}{Im(A)}\sseteq R.
\eeq
 This ideal is $G_{lr}$-invariant and contains the determinantal ideal $I_m(A)$.

\

  The annihilator-of-cokernel is a rather delicate invariant,
   for example it does not transform nicely under  ring homomorphisms, $ann.coker(\phi(A))\neq \phi(ann.coker(A))$. Yet, this ideal is controlled by
 the determinantal ideals as follows:
 \bproperty
\label{Thm.Background.ann.coker.in.terms.of.Fitting.ideals}
\bee
 \item $ann.coker(A)^m\sseteq I_m(A)\sseteq ann.coker(A)\sseteq\sqrt{I_m(A)}$,

  and for any   $m\ge j>1$ holds: $ann.coker(A)\cdot I_{j-1}(A)\sseteq I_j(A)$,
\text{ \cite[proposition 20.7]{Eisenbud}}.
\item  If $m=n$ and $det(A)$ is not a zero divisor in $R$ then $ann.coker(A)=I_m(A):I_{m-1}(A)$.
\item  If $m<n$ and $height(I_m(A))=(n-m+1)$ then $ann.coker(A)=I_m(A)$, \cite[exercise 20.6]{Eisenbud}.
\item If $I_{m-1}(A)$ contains a non-zero divisor modulo $I_m(A)$ then $ann.coker(A)=I_m(A)$, \cite[Corollary 1.4]{Buchsbaum-Eisenbud}.
\eee
\eproperty
In particular, for one-row matrices, $m=1$, or when $I_m(A)$ is a radical ideal, $I_m(A)=ann.coker(A)$.

We use the following properties:
\bel\label{Thm.Background.Ann.Coker.Properties}
\bee
\item Block-diagonal case, $ann.coker(A\oplus B)=ann.coker(A)\cap ann.coker(B)$.
\item For a square matrix,  $m=n$, whose determinant is a non-zero divisor, holds: $ann.coker(A)=ann.coker(A^t)$.

 Suppose, moreover, either $R$ is a unique factorization domain (UFD) or  $height(I_m(A))< height(I_{m-1}(A))$,
 then $ann.coker(A)$ is a principal ideal.
\item If $m>n$ then $ann.coker(A)=\{0\}$.
\eee\eel
The statement {\bf 1} is immediate.

The first part of Statement {\bf 2} follows directly, e.g. from $ann.coker(A)=I_m(A):I_{m-1}(A)$ of property
 (\ref{Thm.Background.ann.coker.in.terms.of.Fitting.ideals}).
For the second part of {\bf 2} note that for a square matrix $I_m(A)$ is a principal ideal and
 the height of $ann.coker(A)\sset R$ is one. If $R$ is UFD then $ann.coker(A)$ is generated by just one element.

Statement {\bf 3}: in this case the submodule $Im(A)\sset R^{\oplus m}$ is of rank$\le n<m$.
 Thus $rank(coker(A))>0$, i.e. this module has linearly independent elements, hence $ann(coker(A))=\{0\}$.
Geometrically, the module $coker(A)$ is supported on the whole $Spec(R)$.

\subsection{Skew-symmetric matrices and their Pfaffian ideals}\label{Sec.Background.Pfaffians}
Suppose $A$ is a square matrix then $I_m(A)=(det(A))$. If moreover $A\in Mat^{skew-sym}_{m\times m}(R)$, with $m$-even,
 then the determinant  is a full square, $det(A)=Pf(A)^2$. Here the Pfaffian polynomial is explicitly written in terms of the entries of $A$,
  \cite{Heymans}, \cite[exercise A.2.11]{Eisenbud}, \cite{Greub}.
In particular, the Pfaffian polynomial transforms under ring homomorphisms, $R\stackrel{\phi}{\to}S$ as $Pf(\phi(A))=\phi(Pf(A))$.

The transformation of the Pfaffian under the congruence is: $Pf(UAU^t)=det(U)Pf(A)$.

If $m$ is odd then $det(A)=0$. In this case, as a measure of degeneracy, we need the following refinement of the Pfaffian ideal:
\beq
(Pf_{m-1}(A)):=\suml_{i=1}^m \big(Pf(A_{\hi})\big).
\eeq
Here $A_{\hi}$ is the $(m-1)\times(m-1)$ block of $A$, obtained by erasing the $i$'th row and column.

We use the following properties of these ideals:
\bel\label{Thm.Background.AntiSymmetric.Matrices.Pfaffians}
Let $R$ be a ring of zero characteristic.
\bee
\item Let $A\in Mat^{skew-sym}_{m\times m}(R)$  for $m$-even. There exists the ``Pfaffian-adjoint" matrix,
 $adj^{Pf}(A)\in Mat^{skew-sym}_{m\times m}(R)$, that satisfies $A\cdot adj^{Pf}(A)=Pf(A)\cdot\one=adj^{Pf}(A)\cdot A$.
 In particular, $ann.coker(A)\supseteq (Pf(A))$.
\item  Moreover: $Span_R(AU+U^tA)_{U\in Mat_{m\times m}(R)}\supseteq Mat^{skew-sym}_{m\times m}(Pf(A))$.
\item  Let $A\in Mat^{skew-sym}_{m\times m}(R)$  for $m>1$, odd. Then
$Span_R(AU+U^tA)_{U\in Mat_{m\times m}(R)}\supseteq Mat^{skew-sym}_{m\times m}(Pf_{m-1}(A))$.
\eee
\eel
\bpr
{\bf 1.} We recall the relevant exterior algebra, \cite[exercise A.2.11]{Eisenbud}. Let $V$ be a free $R$-module of rank $m$.
 Its dual, $V^*$, is the module of one-forms, take the dual basis $\{dx_i\}$.
Associate to $A$ the skew-symmetric form $w_A:=\suml_{1\le i<j\le m}a_{ij}dx_i\wedge dx_j\in\overset{2}{\wedge} V^*$.
 The top exterior power of $w_A$ gives the Pfaffian:
\beq
\overset{\frac{m}{2}}{\wedge}w_A=m!\cdot Pf(A)dx_1\wedge\dots\wedge dx_m\in\overset{m}{\wedge} V^*.
 \eeq
Now, using the pairing  $\overset{2}{\wedge} V^*\otimes \overset{m-2}{\wedge} V^*\to \overset{m}{\wedge} V^*$, we associate to the
  form $\frac{1}{m!}\overset{\frac{m}{2}-1}{\wedge}w_A$ an element of $\overset{2}{\wedge} V^*$, which we call $adj^{Pf}(A)$.
   Then the relation $w_A\wedge (\overset{\frac{m}{2}-1}{\wedge}w_A)=Pf(A)dx_1\wedge\dots\wedge dx_m$ transforms into $A\cdot adj^{Pf}(A)=Pf(A)\cdot\one$.

{\bf 2.} This is now immediate, one takes $U$ in the form $adj^{Pf}(A)\cdot\tilde{U}$, where $\tilde{U}\in Mat^{skew-sym}_{m\times m}(R)$.

{\bf 3.} As in part 1 we have:
\beq
\overset{\frac{m-1}{2}}{\wedge}w_A=
(m-1)!\suml^m_{i=1}Pf(A_\hi)dx_1\wedge\dots\wedge\underbrace{(dx_i)}_{omitted}\wedge\dots\wedge dx_m\in\overset{m-1}{\wedge} V^*.
 \eeq
We get: $dx_i\wedge\Big(\overset{\frac{m-3}{2}}{\wedge}w_A\Big)\wedge w_A=(m-1)!Pf(A_\hi)dx_1\wedge\dots\wedge dx_m$.


Now, we use the pairing  $\overset{2}{\wedge} V^*\otimes \overset{m-3}{\wedge} V^*\to \overset{m-1}{\wedge} V^*$, and the collection of
 forms $\{(\overset{\frac{m-3}{2}}{\wedge}w_A)\wedge dx_i\}_i$ to get:
\beq
Span_R(U_iA+AU_i^t)_i\supseteq (Pf_{m-1}(A))\cdot \overset{2}{\wedge} V^*=Mat^{skew-sym}_{m\times m}(Pf_{m-1}(A)).
\eeq
(Here the matrix $U_i$ corresponds to $\overset{\frac{m-1}{2}}{\wedge}w_A\wedge dx_i$.)
\epr

\subsection{Integral closure of ideals}\label{Sec.Background.Integral.Closure.Ideals}
The integral closure of an ideal $I\sset R$ is defined as
\beq
\bI=\Big\{f\in R\big| \ f^d+\suml^d_{j=1}a_j f^{d-j}=0,\ \text{for some }d\in\N \text{ and some } \{a_j\in I^j\}_{j=1,\dots,d}\Big\}.
\eeq
 This subset is an ideal in $R$ and is itself integrally closed, i.e. $\bar{\bI}=\bI$.
We mention several useful properties:
\bproperty\label{Thm.Background.Integral.Closure.Properties}
\bee
\item $I\sseteq\overline{I}\sseteq\sqrt{I}$. \quad If $I\sset R$ is radical then $\overline{I}=I$.

 If the height of $I\sset R$ is defined then $height(I)=height(\overline{I})$;
\item If $R$ is normal (i.e. integrally closed in its total ring of fractions) then every principal ideal is integrally closed,  \cite[Proposition 1.5.2]{Huneke-Swanson}.
\item $\overline{I:J}\sseteq\overline{I}:J$, \cite[page 7]{Huneke-Swanson}.
\item If $R$ is Noetherian then $\overline{I}\cdot \overline{J}\sseteq \overline{I\cdot J}$, \cite[Corollary 6.8.6]{Huneke-Swanson}.
\item If $R$ is a domain and $\{0\} \neq I$ is finitely generated then $\overline{I\cdot J}:\overline{I}=\overline{J}$, \cite[Corollary 6.8.7]{Huneke-Swanson}.
\item If $R$ is local, Noetherian and regular then $\overline{f\cdot J}=(f)\cdot \overline{J}$.
\item Let $(R,\cm)$ local Noetherian and $J\sset \cm$ an $\cm$-primary ideal. Then $\overline{I}\cdot\hR=\overline{I\cdot \hR}$,
\cite[Lemma  9.3.1]{Huneke-Swanson}.
\eee
\eproperty

In general the computation of $\bI$ is algorithmically involved, it is realized in computer algebra packages as e.g. SINGULAR.
The following  criterion is highly useful in many cases, \cite[Theorem 6.8.3]{Huneke-Swanson}.
 For any ring homomorphism $R\stackrel{\phi}{\to}S$ denote by $S\phi(I)$ the ideal generated by the image of $I$ in $S$.
\beq\label{Eq.Integr.Closur.Val.Criter}
\text{If $R$ is Noetherian then }
\bI=\bigcapl_{R\stackrel{\phi}{\to}S_{DVR}}\phi^{-1} S\phi(I)
\eeq
(Here the intersection goes over all the homomorphisms to  discrete valuation rings.)
In words:
\beq
\text{$f\in \overline{I}$ iff for any homomorphism to a DVR, $R\stackrel{\phi}{\to}S,\ \phi(f)\in S\phi(I)$.}
\eeq
This criterion has a simple geometric formulation (initially  for $R=\C\{\ux\}$, \cite[1.3.4]{Teissier}):
\beq
f\in \overline{I} \text{ iff for any map of the smooth curve-germ, } (C,0)\stackrel{\nu}{\to}Spec(R),
\text{ the pullbacks satisfy: } \nu^*(f)\in\nu^*(I).
\eeq

\subsection{Integral closure of modules}\label{Sec.Background.Integral.Closure.Modules}

For the general definition of the integral closure of modules see \cite[definition 16.1.1]{Huneke-Swanson}. In our case
 all the modules are embedded into $R^{\oplus m}$, the ring is Noetherian over a field with $char(\k)=0$, thus the definition simplifies.
  For any ring homomorphism $R\stackrel{\phi}{\to}S$
 take the image $\phi(M)\sseteq S^{\oplus m}$ and denote by $S\phi(M)$ the generated submodule.
\bed\cite[\S16.1]{Huneke-Swanson}
The integral closure of a submodule $M\sset\R^m$ is
$\overline{M}:=\bigcapl_{R\stackrel{\phi}{\to}S_{DVR}}\phi^{-1} S\phi(M)$.
\eed
In words: $z\in\overline{M}$ iff for any homomorphism to a discrete valuation domain, $R\stackrel{\phi}{\to}S_{DVR}$, there holds   $\phi(z)\in S\phi(M)$.

\

For any ideals $J_1,J_2\sset R$ holds,  \cite[Remark 1.3.2]{Huneke-Swanson}:
\beq\label{Eq.Background.integral.closure.product.of.ideals}
J_1\cdot J_2 \sseteq J_1\cdot\overline{J_2}\sseteq \overline{J_1}\cdot\overline{J_2}\sseteq \overline{J_1\cdot J_2}.
\eeq
More generally:

\bel\label{Thm.Background.Integral.Closure.Lemma}
Let $R$ be  Noetherian, $N\sseteq M\sseteq R^{\oplus m}$ embedded modules, $J\sset R$ an ideal.
\bee
\item $J\cdot\overline{M}\sseteq \overline{J}\cdot\overline{M}\sseteq \overline{J\cdot M}$.
\item $ann\quots{\overline{M}}{\overline{N}}\supseteq\overline{ann\quots{M}{N}}$.
\eee
\eel
\bpr
1. The first inclusion is trivial. For the second it is enough to prove: for any $f\in \overline{J}$, $z\in\overline{M}$ holds $f\cdot z\in \overline{J\cdot M}$.
 Let $R\stackrel{\phi}{\to}S_{DVR}$ then $\phi(f)\in S\cdot \phi(J)$ and $\phi(z)\in S\cdot\phi(M)$.
  Thus $\phi(fz)=\phi(f)\phi(z)\in S\cdot \phi(J)\cdot\phi(M)=S\cdot \phi(J\cdot M)$. As this holds for any such $\phi$ we get $f\cdot z\in \overline{J\cdot M}$.

2. Let $f\in \overline{ann\quots{M}{N}}$, i.e. for any
$R\stackrel{\phi}{\to}S_{DVR}$ holds $\phi(f)\in
\phi(ann\quots{M}{N})$. Thus $f\in ann\quots{M}{N}+ker(\phi)$, i.e.
$f\cdot M\sseteq N+ker(\phi)$.
 But then $\phi(f)\phi(M)\sseteq\phi(N)$. As this holds for any $\phi$, we get: $\overline{f\cdot M}\sseteq\overline{N}$.
  Finally, by part (1.), $f\cdot \overline{M}\sseteq\overline{f\cdot M}\sseteq\overline{N}$, i.e. $f\in ann\quots{\overline{M}}{\overline{N}}$.
\epr
\bex
The ideals in part (2) of Lemma \ref{Thm.Background.Integral.Closure.Lemma} can differ significantly. For example, let $R=\k[[x_1,\dots,x_p]]$,
 with $p\ge3$, and let $N=(x^d_1,\dots,x^d_p)\subsetneq M=\cm^d$.
 Then $\overline{N}=\overline{M}=\cm^d$, hence $ann\quots{\overline{M}}{\overline{N}}=R$. But $\overline{ann\quots{M}{N}}=\cm^{d}$.
\\Proof: obviously $N\sseteq ann\quots{M}{N}$,
 thus $\overline{ann\quots{M}{N}}\supseteq\overline{N}=\cm^{d}$. Furthermore, $\overline{ann\quots{M}{N}}$ is a monomial ideal.
 Suppose $i_1+\cdots +i_p=d-1$ and
$x^{i_1}_1\cdots x^{i_p}_p\in \overline{ann\quots{M}{N}}\smin \cm^d$.
Take the monomial
 $x^{d-1-i_1}\cdots x^{d-1-i_p}_p\in
\cm^{p(d-1)-(d-1)}\sset\cm^d=M$.
 Then
 \[
 (x^{d-1-i_1}\cdots x^{d-1-i_p}_p)\cdot (x^{i_1}_1\cdots x^{i_p}_p)=x^{d-1}_1\cdots x^{d-1}_p\not\in N.
 \]
 Hence $(x^{i_1}_1\cdots x^{i_p}_p)\cdot M\not\subseteq N$, contradicting the assumption  $x^{i_1}_1\cdots x^{i_p}_p\in \overline{ann\quots{M}{N}}$.
\eex

The following lemma compares the bounds in part three of Theorem \ref{Thm.Results.Annihilators.for.GlGrGlr}.
\bel\label{Thm.Background.Comparing.ann.coker.vs.integral.closures}
If $R$ is a Noetherian ring then  $\overline{ann.coker(A)}\sseteq \overline{I_m(A)}:\overline{I_{m-1}(A)}$.
\eel
\noindent Indeed, by Property \eqref{Thm.Background.ann.coker.in.terms.of.Fitting.ideals} we have $ann.coker(A)\cdot I_{m-1}(A)\sseteq I_m(A)$.
 Thus the statement follows from part 1 of Lemma \ref{Thm.Background.Integral.Closure.Lemma}.

\subsection{The largest ``reasonable" group is $G_{lr}$ (matrix preservers over rings)}\label{Sec.Background.Largest.Reasonable.Group}
When looking for the possible groups acting on matrices, $G\circlearrowright\Mat$,
 the guiding principle is that the action of $G$ should, at least, distinguish between degenerate and
 non-degenerate matrices. Groups with such properties are restricted by the following proposition.
 For the general introduction to the theory of preservers, i.e. self-maps of $Mat_{m\times n}(\k)$, that preserve
 some properties/structures see \cite{Molnar}, \cite{Li-Pierce}.
\bprop
Let $R$ be a unique factorization domain with at least four elements. Let $T$ be an invertible map acting $R$-linearly on the module  $Mat_{m\times n}(R)$, for $m\le n$, possibly violating the row/column structure.
Suppose $T$ preserves the set of degenerate matrices, i.e.  $I_m(A)=\{0\}$ iff $I_m(T(A))=\{0\}$.
Then either $T(A)=UAV$, for some $U\in GL(m,R)$, $V\in GL(n,R)$, or, when $m=n$, $T(A)=UA^tV$.
\eprop
\bpr
We pass to the field of fractions $Frac(R)$. As $R$ has no zero-divisors, it is naturally embedded, $R\sset Frac(R)$.
Thus $T$ extends to an invertible operator that acts $Frac(R)$-linearly on
the vector space of matrices $Mat_{m\times n}(Frac(R))$
 and such that $I_m(A)=\{0\}$ if and only if $I_m(T(A))=\{0\}$.
 In other words, $T$ preserves the degenerate matrices, of $rank<m$, and also preserves the non-degenerate matrices  of $rank=m$.
Now we use the classical results on the preservers over a field:
\bee[i.]
\item(the case $m=n$) \cite[Theorem 3]{Dieudonne}: Let $\k$ be an arbitrary field. Any invertible $\k$-linear transformation of $Mat_{m\times m}(\k)$ that preserves the set of degenerate matrices has either the form   $T(A)=UAV$ or the form  $T(A)=UA^tV$, for some $U,V\in GL(m,\k)$.
\item(the case $m<n$) \cite[Theorem 7]{Beasley-Laffey} Suppose the field $\k$ contains at least four elements. Any invertible $\k$-linear transformation of $Mat_{m\times n}(\k)$ that preserves the set of $rank=r$ matrices (for some fixed $0<r\le m$) is of the form $T(A)=UAV$  for some $U\in GL(m,\k)$, $V\in GL(n,\k)$.
\eee

Therefore we get: either $T(A)=UAV$, for some $U\in GL(m,Frac(R))$, $V\in GL(m,Frac(R))$, or, for $m=n$, $T(A)=UA^tV$.
 It remains to show that we can choose $U\in GL(m,R)$, $V\in GL(n,R)$. (We check the case $m<n$, the case $m=n$ is similar.)

Let $0\neq g\in R$  be a(ny) common denominator of the entries of $U$, i.e. $g\cdot U\in Mat_{m\times n}(R)$. Replace $(U,V)$ by $(g\cdot U,\frac{V}{g})$,
 so we still have $T(A)=UAV$.

 Substitute for $A$  the elementary matrix $E_{ij}$, with 1 at the i'th row, j'th column and zeros otherwise.
 The $p,q$-entry of the matrix $UE_{ij}V$ equals $u_{pi}v_{jq}$ and it belongs to $R$. Suppose some element of $V$ is not in $R$, we can assume
 it is $v_{11}$. Present $v_{11}$ as a fraction,  $v_{11}=\frac{f}{g}$, for some $f,g\in R$. As $R$ is UFD we can assume that $f,g$ are co-prime,
  i.e. $(f)\cap (g)=(fg)$ and $g$ is not invertible in $R$. As $u_{pi}v_{11}\in R$, for any $p,i$, we get: $U$ is divisible by $g$.
  Thus we pass from $(U,V)$  to $(\frac{U}{g},gV)$ and proceed to the other entries of $V$.
 \epr

Therefore, in our work the group $G\circlearrowright\Mat$ is always a subgroup of $G_{lr}=GL_R(m)\times GL_R(n)$.

\beR  The last statement is very sensitive to the particular form of condition ``$T$ preserves degenerate matrices".
 For example, when the ring is local, one could ask for a condition: $T(A)\in Mat_{m\times n}(\cm^q)$ iff $A\in Mat_{m\times n}(\cm^q)$.
 But this condition is satisfied just by $R$-linearity of $T\circlearrowright\Mat$, without any additional restrictions.
\eeR

\subsection{Unipotent subgroups}\label{Sec.Background.Unipotent.Subgroups}
Fix a group action on an $R$-module $G\circlearrowright M$ and a decreasing filtration $M=M_0\supsetneq M_1\supsetneq\cdots$. We get the induced filtration of $G$ by the
 (normal) subgroups
\beq
G^{(i)}:=\{g\in G|\ \forall\ j\ge0\quad \forall\ z_j\in M_j:\ g(z_j)-z_j\in M_{j+i}\}.
\eeq
In particular, the subgroup $G^{(0)}\sseteq G$ preserves the filtration, while the subgroup $G^{(1)}\sseteq G$ acts unipotently. We call $G^{(1)}$ the ``unipotent" subgroup of $G$. If $G^{(1)}=G$ then $G$ itself is called ``unipotent for the given filtration".
 The general properties of the filtration $G\supseteq G^{(0)}\supseteq G^{(1)}\supseteq\cdots$ are studied in \cite[\S3]{Belitski-Kerner}.

 For the goals of the current paper we usually need the filtration by the powers of a/the maximal ideal,
\beq
\Mat\supsetneq \Matm\supsetneq Mat_{m\times n}(\cm^2)\cdots.
\eeq
Therefore we denote the unipotent subgroup by $G^{(\cm)}\sseteq G$. In particular, when $R$ is a local ring:
\beq
G^{(\cm)}_r=\{U=\one+\tilde{U}|\ \tilde{U}\in Mat_{n\times n}(\cm)\},\quad
G^{(\cm)}_l=\{U=\one+\tilde{U}|\ \tilde{U}\in Mat_{m\times m}(\cm)\},\quad
G^{(\cm)}_{lr}=G^{(\cm)}_l\times G^{(\cm)}_r.
\eeq
The last equality is a statement, but the proof is straightforward, the conditions
 $\{UE_{ij}V-E_{ij}\in Mat_{m\times n}(\cm)\}_{i,j}$ force $U\in G^{(\cm)}_l$, $V\in G^{(\cm)}_r$.

More generally, for any subgroup $G\sset G_{lr}$ we have: $G^{(m)}=G\cap G^{(\cm)}_{lr}$.

\subsection{The tangent spaces}\label{Sec.Background.Tangent.Spaces}
The tangent spaces to the groups and group-orbits are defined and studied in great generality in \cite[\S 3.7]{Belitski-Kerner}.
Here we recall them for our group actions.

We use the isomorphism of $R$-modules, $T_{(\Mat,A)}\isom\Mat$, to identify $T_{(GA,A)}$,  $T_{(\Si,A)}$ with their images in $\Mat$.
\bex\label{Ex.Tangent.Spaces.to.the.Actions}  (See \cite[\S 3.8]{Belitski-Kerner}.)
\bee[i.]
\item
$\Si=\Mat$ and $G_{lr}:\ A\to UAV^{-1}$.
Here $T_{(G_{lr}A,A)}=Span_R\{uA,Av\}_{(u,v)\in Mat_{m\times m}(R)\times Mat_{n\times n}(R)}$. Similarly
 for $G_l$ and $G_r$.
\item
$\Si=Mat^{sym}_{m\times m}(R)$ or $Mat^{skew-sym}_{m\times m}(R)$  and $G_{congr}:\ A\to UAU^t$.
 Here $T_{(G_{congr}A,A)}=Span_R\{uA+Au^t\}_{u\in Mat_{m\times m}(R)}$.
\item
$\Si=Mat_{m\times m}(R)$  and $G_{conj}:\ A\to UAU^{-1}$. Here $T_{(G_{conj}A,A)}=Span_R\{uA-Au\}_{u\in Mat_{m\times m}(R)}$.
\item
$\Si=Mat^{up}_{\{m_i\}\times\{n_j\}}(R)$  and $G^{up}_{lr}:\ A\to UAV$, where $U,V$ are upper-block-triangular. Here
\[T_{(G_{lr}A,A)}=Span_R\{uA,Av\}_{(u,v)\in Mat^{up}_{\{m_i\}\times \{m_j\}}(R)\times Mat^{up}_{\{n_i\}\times \{n_j\}}(R)}.\]
 Similarly for  $G^{up}_r$.
\eee
\eex
In all these cases  holds: $T_{(G^{(\cm)}A,A)}=\cm\cdot T_{(GA,A)}$.
In general, for all the ``reasonable" group actions there holds
\beq
\cm\cdot T_{(GA,A)}\subseteq T_{(G^{(\cm)}A,A)}\subseteq T_{(GA,A)},
\eeq
 see \cite[\S3.4]{Belitski-Kerner} for more detail.

\subsection{Transition to the completion}
Fix some proper ideal $\{0\}\neq \ca\ssetneq R$. Here $\ca$ is not necessarily a maximal ideal, as we need these results
 for partial algebraization, \S\ref{Sec.Results.Algebraization},  and relative determinacy, \S\ref{Sec.Results.Relative.Determinacy}.
 We often take the $\ca$-adic completion, $R\stackrel{\hat{\psi}}{\to}\hR$,
 and compute over $\hR$. To translate the results to $R$ we need the following standard facts.

Given two embedded modules, $M_1\sset M_2\sset R^{\oplus p}$, take their images under the completion map,
 $\hat{\psi}(M_1)\sseteq\hat{\psi}(M_2)\sseteq \hat\psi(R)^{\oplus p}$.
  These are not necessarily $\hR$-modules (unless the completion map is surjective), but they generate the $\hR$-modules:
 $\hR\cdot\hat{\psi}(M_1)\sseteq\hR\cdot\hat{\psi}(M_2)\sseteq \hR^{\oplus p}$.

Note that in general $\hat{\psi}(M_i)$, $\hR\cdot\hat{\psi}(M_i)$ differ from the ``non-embedded" completion, $\widehat{M}=\liml_{\leftarrow}\quots{M}{\ca^n\cdot M}$.
(For example, fix a non-zero divisor $f\in\ca^\infty$ and consider $J=(f)$. Then $\hat{\psi}(J)=\{0\}\sset\hR$,
 while $\widehat{J}\approx\hR$, as $\hR$-modules.)

In the statements below we use the Jacobson ideal of the ring, $J(R)$. Recall, for a local ring $J(R)=\cm$.

First we address the behaviour of the Loewy length of an ideal, $ll_R(J)=min\{n|\ J\supseteq\ca^n\}$,  under completion.
\bel
Fix an ideal in a commutative unital ring, $\{0\}\neq \ca\ssetneq R$. Fix some other ideal, $\ca^\infty\sseteq J\sset R$. Suppose:\vspace{-0.2cm}

\bee[i.]
\item either the completion map is surjective, $R\stackrel{\hat{\psi}}{\twoheadrightarrow}\hR$;
\item or both $\ca$ and $J$ are finitely generated (as $R$-modules) and $\ca\sseteq J(R)$.
\eee
Then $J\supseteq\ca^n$ iff $\hR\cdot\hat\psi(J)\supseteq \hR\cdot\hat\psi(\ca)^n$. Equivalently: $ll_R(J)=ll_\hR(\hR\cdot \hat\psi(J))$.
\eel
\bpr
The part $\Rrightarrow$ is immediate, as the ideal $\hR\cdot\hat\psi(\ca)^n$ is generated by the images of generators of $\ca^n$.

$\Lleftarrow$  Suppose $R\stackrel{\hat{\psi}}{\twoheadrightarrow}\hR$, then $\hR\cdot\hat\psi(J)=\hat\psi(J)$.  If
 $\hat\psi(J)\supseteq\hat\psi(\ca)^n$,
 then $J+\ca^\infty\supseteq\ca^n$. (By taking some preimages of the elements of $\hat\psi(\ca)^n$ in $J$.) Thus, as $J\supseteq\ca^\infty$, we get $J\supseteq\ca^n$.

Suppose $\ca,J$ are finitely generated, then the inclusion follows by the standard Nakayama arguments. The completion morphism induces
  $\quots{R}{\ca^N}\isom\quots{\hR}{\hR\cdot\hat\psi(\ca)^N}$, for any $N$, and similarly
   $\quots{J+\ca^N}{\ca^N}\isom\quots{\hR\cdot\hat\psi(J)+\hR\cdot\hat\psi(\ca)^N}{\hR\cdot\hat\psi(a)^N}$.
 Thus  $\hR\cdot\hat\psi(J)\supseteq\hR\cdot\hat\psi(\ca)^n$ implies  $\quots{J+\ca^N}{\ca^N}\supseteq \quots{\ca^n}{\ca^N}$, or, equivalently,
    $J+\ca^N\supseteq\ca^n$, for any $N$.
   Therefore, for  the embedding of (finitely generated) modules $J\sseteq J+\ca^n$, we have $J+\ca(J+\ca^n)\supseteq J+\ca^n$.  And by Nakayama
    we have $J=J+\ca^n$, hence $J\supseteq\ca^n$.
\epr

Now we check the change of the annihilator of the quotient, $ann\quots{M_2}{M_1}$, under the completion.
\bel\label{Thm.Background.Completion.Lemma}
Fix an ideal in a commutative unital ring, $\{0\}\neq \ca\ssetneq R$, and two modules, $M_1\sset M_2\sseteq R^{\oplus p}$,
 with $M_1\supseteq M_2\cap(\ca^\infty\cdot R^{\oplus p})$. Suppose:\vspace{-0.2cm}

\bee[i.]
\item either the completion map is surjective, $R\stackrel{\hat{\psi}}{\twoheadrightarrow}\hR$;
\item or $R$ is Noetherian and $M_2$ is a finitely generated $R$-module and $\ca\sseteq J(R)$.
\eee
Then $ann\quots{M_2}{M_1}\supseteq \ca^n$  iff
 $ann\big(\quots{\hR\cdot\hat{\psi}(M_2)}{\hR\cdot\hat{\psi}(M_1)}\big)\supseteq\hR\cdot\hat{\psi}(\ca^n)$. Equivalently:
    $ll_R(ann\quots{M_2}{M_1})=ll_\hR(ann\quots{\hR\cdot\hat{\psi}(M_2)}{\hR\cdot\hat{\psi}(M_1)})$.
\eel
\bpr
$\Rrightarrow$ Direct check: $\hR\cdot\hat{\psi}(\ca^n)\cdot \hR\cdot\hat{\psi}(M_2)=\hR\cdot\hat{\psi}(\ca^n \cdot M_2)\sseteq \hR\cdot\hat{\psi}(M_1)$.

$\Lleftarrow$ Suppose $\hR\cdot\hat{\psi}(\ca^n)\cdot \hR\cdot\hat{\psi}(M_2)\sseteq \hR\cdot\hat{\psi}(\ca^n \cdot M_1)$.

In the case of surjective completion, $R\stackrel{\hat{\psi}}{\twoheadrightarrow}\hR$, we get $\ca^n M_2\sseteq M_1+\ca^\infty R^{\oplus p}$, by choosing some preimages
 under $\hat\psi$.
 Then, by direct check,
 $\ca^n M_2\sseteq M_1+M_2\cap \ca^\infty R^{\oplus p}=M_1$.

Otherwise, when $R$ is Noetherian, $\ca,M_2$ are finitely generated, and we use the Nakayama argument. As in the previous lemma, we have: $\ca^nM_2\sseteq M_1+(M_2\cap \ca^NR^{\oplus p})$,
 for any $N$. By Artin-Rees lemma, for $N\gg1$ we have: $(M_2\cap \ca^NR^{\oplus p})\sseteq \ca^{n+1}M_2$. Thus we have
  $\ca^nM_2\sseteq M_1+\ca^{n+1}M_2$. Now the Nakayama lemma gives $\ca^nM_2\sseteq M_1$.
\epr

\bel\label{Thm.Background.Completion.Lemma.surjective.case}
As before, fix two modules, $M_1\sset M_2\sseteq R^{\oplus p}$,  with $M_1\supseteq M_2\cap(\ca^\infty\cdot R^{\oplus p})$.
 Suppose the $\ca$-adic completion is surjective, $R\stackrel{\hat{\psi}}{\twoheadrightarrow}\hR$, then
 $\hat{\psi}(ann\quots{M_2}{M_1})=ann\quots{\hat{\psi}(M_2)}{\hat{\psi}(M_1)}$.
\eel
\bpr
The part $\sseteq$ is straightforward: $\hat\psi(ann\quots{M_2}{M_1})\cdot\hat\psi(M_2)=\hat\psi(ann(\quots{M_2}{M_1})\cdot M_2)\sseteq \hat\psi(M_1)$.

The part $\supseteq$.
By the surjectivity $R\twoheadrightarrow\hR$  any element of $ann\quots{\hat{\psi}(M_2)}{\hat{\psi}(M_1)}$ is of the form $\hat{\psi}(f)$ for some $f\in R$. As $\hat{\psi}(f)\hat{\psi}(M_2)\sseteq\hat{\psi}(M_1)$ we get: $fM_2\sseteq M_1+(M_2\cap\ca^\infty R^{\oplus p})=M_1$. Thus $f\in ann\quots{M_2}{M_1}$.
\epr

\bcor\label{Thm.Loewy.Length.under.completion}
Fix two ideals, $I,J\sset R$, with $I\supseteq\Big((I+J)\cap\ca^\infty\Big)$.
\bee
\item
Suppose the completion is surjective, $R\twoheadrightarrow\hR$,
then $\hat{\psi}(I:J)=\hat{\psi}(I):\hat{\psi}(J)$.

\item Suppose $R$ is Noetherian and $\ca\sseteq J(R)$, then  $ll_R(I:J)=ll_\hR(\hR\hat\psi(I): \hR\hat\psi(j))$.

\item Suppose $(R,\cm)$ is a local Noetherian ring and $A\in \Mat$ satisfies: $I_m(A)$ is an $\cm$-primary ideal. Then, for the $\cm$-completion $R\stackrel{\hat\psi}{\to}\hR$
 holds: $ll_R\Big(\overline{I_m(A)}:\overline{I_{m-1}(A)}\Big)=ll_\hR\Big(\overline{I_m(\hA)}:\overline{I_{m-1}(\hA)}\Big)$.
\eee\ecor
\bpr
Note that $I:J=ann\quots{I+J}{I}$.

1.
Apply part one of Lemma \ref{Thm.Background.Completion.Lemma.surjective.case} to the modules $I\sseteq I+J\sset R$.  We get:
$\hat{\psi}\Big(ann\quots{I+J}{I}\Big)=ann\Big(\quots{\hat{\psi}(I+J)}{\hat{\psi}(I)}\Big)$, hence the statement.

2. Follows by Lemma \ref{Thm.Background.Completion.Lemma}.

3.
By part 1 we have: $ll_R\Big(\overline{I_m(A)}:\overline{I_{m-1}(A)}\Big)= ll_\hR\Big((\overline{I_m(A)}\cdot\hR):(\overline{I_{m-1}(A)}\cdot\hR)\Big)$.
 By the integral closure properties, Property \ref{Thm.Background.Integral.Closure.Properties}, we have:
  $\overline{I_j(A)}\cdot\hR=\overline{I_j(A)\cdot\hR}$. Finally, $I_j(A)\cdot\hR=I_j(A\cdot\hR)=I_j(\hA)$.
\epr

\subsection{The ``relevant approximation" property}\label{Sec.Background.Approximation.Properties}
 To establish finite determinacy we should resolve the condition $gz=z+w$, $g\in G$. In most cases this is a (finite) system of equations.
  First one looks for an order-by-order solution (modulo the powers of the maximal ideal). Then one tries to establish an ordinary solution for these equations.

For the criteria in the Subsection \ref{Sec.Background.BK1} we need the ``relevant approximation property of the ring $R$". The particular property depends on
 the type of the equations.
While the general theory is developed in \cite{Belitski-Kerner},
for the subgroups of $G_{lr}$ we need only the following two types of approximation property.

\subsubsection{} If $G\sseteq G_{lr}$ is defined by $R$-linear equations  and the condition $gA=A+B$ can be written as a system of linear equations on $g=(U,V)$ then
we say that $R$ has the relevant approximation property if
 the assumption \eqref{Eq.Intro.assumption.R.Noetherian.or.surjects} together with the condition $ann.coker(A)\supseteq\cm^\infty$ hold.

\

The equations are linear for the groups $G_r$, $G_l$, $G_{lr}$, $C_{conj}$.
They are:
\beq
UA=A+B,\quad\quad AV=A+B,\quad\quad UA=(A+B)V,\quad\quad UA=(A+B)U.
\eeq
The class of rings satisfying \eqref{Eq.Intro.assumption.R.Noetherian.or.surjects} is rather large. Besides the traditional Noetherian
 rings of Commutative Algebra this class contains e.g.  the ring $C^\infty(\R^p,0)$ of infinitely differentiable function-germs in several variables.
   Its completion is
  Noetherian, $\R[\![\ux]\!]$, and the completion map, $C^\infty(\R^p,0)\twoheadrightarrow\R[\![\ux]\!]$, is surjective
  by Borel's lemma, \cite[pg. 284, exercise 12]{Rudin-book}.

For the non-Noetherian rings, with $\cm^\infty\neq\{0\}$ the condition $ann.coker(A)\supseteq\cm^\infty$ can be checked using the classical
{\L}ojasiewicz property:
\begin{multline}\label{Eq.Lojasiewicz.Ineq}
\text{If $f\in C^\infty(\R^p,0)$ satisfies $|f|\ge C|\ux|^\de$, for some $C>0$, $\de>0$,}
 \\\text{then $f$ divides any function flat at the origin, and thus $(f)\supseteq\cm^\infty$.}
\end{multline}

\subsubsection{} If $G\sseteq G_{lr}$ is defined by polynomial/analytic equations and the condition $gA=A+B$ can be written as a system of polynomial/analytic
equations on $g=(U,V)$,
 we say that $R$ has the relevant approximation property provided $R$ is Henselian ring.
 This is the case for $G_{congr}$, the equations being quadratic, $UAU^t=A+B$.

\subsection{The transition from $ann(T^1_{(\Si,G,A)})$ to the finite determinacy and the main results of \cite{Belitski-Kerner}}\label{Sec.Background.BK1}
In this section $(R,\cm)$ is a local ring over a field, $\k$, of zero characteristic.

Fix a finitely generated $R$-module with descending filtration, $M=M_0\supset M_1\supset\cdots$.
 We assume that the filtration is ``essentially descending", i.e. for any $j$ there exists $k_j$ such that
   $M_{k_j}\sseteq \cm^j\cdot M$.
Fix a (filtered) group action, $G\circlearrowright M$, and a deformation space $\Si\sseteq M$. (In our case $M=\Mat$ or $Mat^{sym}_{m\times m}(R)$
 or  $Mat^{skew-sym}_{m\times m}(R)$, and $\Si\sseteq M$ is a submodule.)
 To understand the ``essential" deformations of an element $z\in \Si$, we need to know
  ``how small/large" is the orbit $Gz$ as compared to $\Si$. In many cases this question reduces to the linear version:
\beq\label{Eq.question}
\text{What is the biggest $R$-submodule $\La\sset M$ satisfying: $z+\La\sseteq Gz$?}
\eeq
(The precise statement is \cite[Lemma 2.1]{Belitski-Kerner}.)

The main result of \cite{Belitski-Kerner} is the transition of this linear version to the comparison of the
 tangent spaces, $T_{(Gz,z)}\sseteq T_{(M,z)}$. We assume that the germ $(Gz,z)$ is ``nice", in particular
 it possesses  a well defined tangent space. The precise condition is:
\beq\label{Eq.assumptions.of.kpd.fs.etc}\ber
\text{\em the subgroup $G\sset GL_\k(M)$ and its completion $\widehat{G}\sseteq GL_\k(\hM)$ are \kpd;}\\
\text{\em their unipotent parts, $G^{(1)}$, $\widehat{G}^{(1)}$, are of Lie type}.
\eer\eeq
The condition  ``$G$ is of Lie type" ensures that the tangent space $T_{(Gz,z)}$ ``approximates" the germ $(Gz,z)$.
 These properties hold for all our groups, $G_l,G_r,G_{lr}, G_{congr}, G_{conj},\dots$,  see the details in \cite[\S3.7]{Belitski-Kerner}.

Furthermore, we use the isomorphism of $R$-modules $T_{(M,z)}\isom M$ to identify $T_{(Gz,z)}$ with its embedding.
 Accordingly, we have the embedded module: $T_{(Gz,z)}\sseteq M$, as in \S\ref{Sec.Background.Tangent.Spaces}.

\bthe\label{Thm.Background.Linearization}\cite[Theorem 2.2]{Belitski-Kerner}
Suppose the (filtered) action $G\circlearrowright\{M_i\}$ satisfies assumptions \eqref{Eq.assumptions.of.kpd.fs.etc}.
 Suppose that $G$ is unipotent for the filtration $\{M_i\}$.
\bee
\item If $M_i\sseteq T_{(Gz,z)}$ and $R$ has the relevant approximation property (\S\ref{Sec.Background.Approximation.Properties}) then $\{z\}+M_i\sseteq Gz$.
\item Suppose $T_{(Gz,z)}\sseteq T_{(M,z)}$ is an $R$-submodule. If $\{z\}+M_i\sseteq Gz$ then $M_i\sseteq T_{(Gz,z)}$.
\eee
\ethe

While Theorem \ref{Thm.Background.Linearization} is quite general,  to compute/bound the order of determinacy we use a more specific
 criterion (see \S\ref{Sec.Results.Notations} for the Loewy length $ll_R(..)$):
\bprop\label{Thm.Intro.ord.of.det.via.Loewy.length}
Suppose $\Si\sseteq\Mat$ is a free direct summand, i.e. $\Si\oplus\Si^\bot=\Mat$ for a free submodule $\Si^\bot\sset\Mat$.
\bee
\item \cite[Corollary 2.6]{Belitski-Kerner} Suppose $ann(T^1_{(\Si,G,A)})\supseteq\cm^\infty$ and $R$ has the relevant approximation property. Then
\[
ll_R\Big(ann(T^1_{(\Si,G,A)})\Big)-1\le ord^\Si_G(A)\le
ll_R\Big(ann(T^1_{(\Si,G^{(\cm)},A)})\Big)-1.
\]

\item If $ann(T^1_{(\Si,G,A)})\not\supseteq\cm^\infty$, then $A$ is not infinitely-$(\Si,G)$-determined.
\eee
\eprop
For all our groups holds
\beq
\cm\cdot T_{(GA,A)}\sseteq T_{(G^{(\cm)}A,A)}\sseteq T_{(GA,A)},
\eeq
 thus we have the bounds on the annihilator:
\beq
\cm\cdot ann(T^1_{(\Si,G,A)})\sseteq ann(T^1_{(\Si,G^{(\cm)},A)})\sseteq ann(T^1_{(\Si,G,A)}).
\eeq
 So, the bounds in the corollary differ at most by one. Moreover, in many cases $ann(T^1_{(\Si,G,A)})=ann(T^1_{(\Si,G^{(\cm)},A)})$, see \S\ref{Sec.Background.Tangent.Spaces}.

This proposition reduces the determinacy question to the study of $T^1_{(\Si,G,A)}$ and its annihilator.

\beR
We mention the relation to the Castelnuovo-Mumford regularity of a module. Recall, \cite[exercise 20.18]{Eisenbud}:
if $M=\oplusl_i M_i$ is a graded $R$-module of finite length then its regularity equals $reg(M)=max\{i|\ M_i\neq\{0\}\}$.
Therefore, when $T^1_{(\Si,G,A)}$ is graded and of finite length, we can rewrite:
\beq
reg(T^1_{(\Si,G,A)})\le ord^\Si_G(A)\le reg(T^1_{(\Si,G^{(\cm)},A)}).
\eeq
\eeR

\subsection{Finite determinacy vs stable algebraizability}\label{Sec.Background.Fin.Determin.vs.Stable.Algebraizability}
Let $\quots{\k[\ux]}{I}\subseteq R\subsetneq \hR=\quots{\k[[\ux]]}{\k[[\ux]]\cdot I}$ be a ring over the field $\k$ of zero characteristic.
  The typical examples of $R$
 are: $\quots{\k[\ux]}{I}$, $\quots{\k\bl\ux\br}{I}$ (algebraic power series), $\quots{\k\{\ux\}}{I}$ (convergent power series).
 In this section we do not use any matrix structure, thus instead of the modules $\Mat$, $Mat_{m\times n}(\hR)$ we work with just free modules
  $F\sset \hF:=\hR\cdot F$.
 By the assumption on $R$ the completion is injective, thus $F\sset \hF$.

Fix a group action $G\circlearrowright \hF$,  $\k$-linear but not necessarily $\hR$-linear.
 We assume $G$ to be complete \wrt the filtration $\{G^{(i)}\}$ and moreover of Lie-type (see \S\ref{Sec.Background.BK1}).
 Then the unipotent part $G^{(1)}$ admits the surjective exponential map, $T_{(G^{(1)},\one)}\stackrel{exp}{\to}G^{(1)}$,
  $\xi\to exp(\xi)=\suml_{j=0}^\infty\frac{\xi^j}{j!}$, see \cite[\S3.5]{Belitski-Kerner}.
 This map sends an infinitesimal family of tangent vectors, $\xi(t)=\sum \xi_j t^j\in T_{(G^{(1)},\one)}[[t]]$, to infinitesimal family of group elements,
  $exp(\xi(t))$. Denote the collection of the later families by $G[[t]]$.

 \bed
 \bee[i.]
\item  An element $\hv\in \hF$ is called $R$-algebraizable if $\hv\in Gv$ for some $v\in F$.
\item  An element $v\in \hF$ is called stably-$R$-algebraizable if there exists $N\in\N$ such that for any deformation $\hv+\hu(t)$, with
$\hu(t)\in (t)\cdot \cm^N\cdot \hF[[t]]$,
   there exists a family $g(t)\in G[[t]]$ satisfying $g(t)(\hv+\hu(t))\in F[[t]]$, and, moreover, if $\hv\in F$ then one can assume $g(0)=\one$.
\eee
\eed
These algebraizability notions depend on the choice of $R$, the strongest notion is for $R=\quots{\k[\ux]}{I}$.

  Finite determinacy is stronger than $\quots{\k[\ux]}{I}$-algebraizability and is equivalent to the stable $\quots{\k[\ux]}{I}$-algebraizability.
\bprop
Let $G\circlearrowright \hF$ be an action of a Lie-type group, suppose $G$ is complete \wrt the filtration $\{G^{(i)}\}$
  and for any $\hv\in \hF$ the tangent space $T_{(G\hv,\hv)}$ is an $\hR$-module.
   Then $\hv\in\hF$ is $G$-finitely determined iff $\hv$ is stably-$R$-algebraizable.
\eprop
\bpr
$\Rrightarrow$ As $\hv$ is finitely determined, we can cut off its Tailor tail, thus we assume $\hv\in F$, and denote it by $v$.
Choose any $N$ bigger than the $G$-order of  determinacy of $v$ and take any deformation $\hu(t)=\suml_{i\ge 1}\hu_it^i$,
 with $\hu_i\in\cm^N\cdot \hF$.
 Thus $\hu_i\in T_{(G^{(1)}v,v)}$ and we fix some $\widetilde{\xi_i}\in T_{(G^{(1)}v,v)}$ satisfying $\widetilde{\xi_i}(v)=u_i$. Thus
\beq
v+\suml_{i\ge 1}u_it^i=(\one+\suml_{i\ge 1}\widetilde{\xi_i}t^i)v=exp\Big(ln(\one+\suml_{i\ge 1}\widetilde{\xi_i}t^i)\Big)v=
  exp\Big( \suml_{i\ge 1}\xi_it^i\Big)v.
  \eeq
   Here $\{\xi_i\}_i$ are the elements of $T_{(G^{(1)},\one)}$ determined by $\{\widetilde{\xi_i}\}_i$.
 Then $exp\Big( \suml_{i\ge 1}\xi_it^i\Big)\in G[[t]]$   is the needed element.

$\Lleftarrow$  We prove $T_{(G^{(1)}v,v)}\supseteq \cm^N\cdot \hF$, by Theorem \ref{Thm.Background.Linearization}
 this implies the finite determinacy.

  We can assume $\hv\in F$, thus denote it by $v$.
 By the assumption, the family $(v+\hu\cdot t)\in \hF[[t]]$, with $\hu\in \cm^N\cdot\hF$, is presentable
 as $exp\Big( \suml_{i\ge 1}\xi_it^i\Big)v$ for some $\{\xi_i\in T_{(G^{(1)},\one)}\}$.
 Therefore $\xi_1 v=\hu$, i.e. $\hu\in T_{(Gv,v)}$.
Now invoke the finite determinacy statement, Theorem \ref{Thm.Background.Linearization}.
 \epr

Note that in this proof we could not just substitute $t=1$ as the element $exp(\sum \xi_i)$  is not well defined.

\section{Proofs, examples and corollaries}\label{Sec.Proofs}
\subsection{The case of $G_l,G_r,G_{lr}\circlearrowright \Mat$}\label{Sec.Proofs.Glr}
\subsubsection{Approximation of $ann(T^1_{(\Si,G,A)})$} \

{\em Proof of Theorem \ref{Thm.Results.Annihilators.for.GlGrGlr}}

{\bf 1.} Note that $T_{(G_rA,A)}=Span_R(AU)_{U\in Mat_{n\times n}(R)}=\oplusl^n_{i=1} Im(A)=Im(A)^{\oplus n}$, where $Im(A)\sset R^{\oplus m}$.
 In addition, $T_{(\Si,A)}=\oplusl^n_{i=1} R^{\oplus m}$. Finally, the two decompositions are compatible, i.e.
 $\quots{T_{(\Si,A)}}{T_{(G_rA,A)}}=\oplusl^n_{i=1}\quots{R^{\oplus m}}{Im(A)}=\oplusl^n_{i=1} coker(A)$.

{\bf 2.} The $(G_l,A)$ case is identical to the $(G_r,A^t)$ case. If $m<n$ then $ann.coker(A^t)=\{0\}$, by Lemma \ref{Thm.Background.Ann.Coker.Properties}.

 Similarly for $m=n$: $ann(T^1_{(\Si,G_l,A)})=ann(T^1_{(\Si,G_r,A^t)})=ann.coker(A^t)=ann.coker(A)$, again, Lemma \ref{Thm.Background.Ann.Coker.Properties}.

{\bf 3.}
 As $G_{lr}\supset G_r$, we have $T_{(G_{lr}A,A)}\supseteq T_{(G_rA,A)}$, thus part 1 gives the bound
\beq
ann\big(T^1_{(\Si,G_{lr},A)}\big)\supseteq ann\big(T^1_{(\Si,G_r,A)}\big)=ann.coker(A).
\eeq

Suppose $f\in ann\big(T^1_{(\Si,G_{lr},A)}\big)$, i.e.
\beq
f\cdot \Mat\sseteq Span_R\Big(UA+AV|\ (U,V)\in Mat_{m\times m}(R)\times Mat_{n\times n}(R)\Big).
\eeq
 For any projection onto a discrete valuation domain, $R\stackrel{\phi}{\to}S$, we get:
\beq
\phi(f)\cdot Mat_{m\times n}(S)\sseteq Span_S\Big(U\phi(A)+\phi(A)V|\ (U,V)\in Mat_{m\times m}(S)\times Mat_{n\times n}(S)\Big).
\eeq
Now we use the canonical form, Proposition \ref{Thm.Background.Canonical.Form.Matrix.over.DVR}, $\phi(A)=PD_{\phi(A)}Q$, where $D_{\phi(A)}$ has zeros off the main diagonal.
We get:
\beq
P^{-1}\phi(f)\cdot Mat_{m\times n}(S)Q^{-1}\sseteq Span_S\Big(UD_{\phi(A)}+D_{\phi(A)}V|\ (U,V)\in Mat_{m\times m}(S)\times Mat_{n\times n}(S)\Big).
\eeq
On the right hand side one has the space of matrices whose $(m,m)$'th component lies in the ideal $(\la_m)$.
 On the left hand side one has $Mat_{m\times n}(\phi(f))$.
Therefore, for $I_{m-1}(\phi(A))\neq\{0\}$:
\beq
\phi(f)\in (\la_m)=\Big( I_m(\phi(A)):I_{m-1}(\phi(A))\Big)=\Big(\phi\big(I_m(A)\big):\phi\big(I_{m-1}(A)\big)\Big).
\eeq
Equivalently one has:
 $\phi(f\cdot I_{m-1}(A))\sseteq\phi(I_m(A))$.  Note that this embedding holds also when $I_{m-1}(\phi(A))=\{0\}$.
As this holds for any $\phi$, we get: $\overline{f\cdot I_{m-1}(A)}\sseteq\overline{I_m(A)}$. As this holds for any $f\in ann\big(T^1_{(\Si,G_{lr},A)}\big)$,
we have:
\beq
\overline{I_{m-1}(A)\cdot ann\big(T^1_{(\Si,G_{lr},A)}\big)}\sseteq\overline{I_m(A)}.
\eeq
Finally, by Lemma \ref{Thm.Background.Integral.Closure.Lemma}, for $I_{m-1}(A)\neq\{0\}$,
 we have:
$\overline{ann\big(T^1_{(\Si,G_{lr},A)}\big)}\sseteq \overline{I_m(A)}:\overline{I_{m-1}(A)}$.

For $I_{m-1}(A)=\{0\}$ similar arguments show: $\phi(ann\big(T^1_{(\Si,G_{lr},A)}\big))=\{0\}$ for any $R\stackrel{\phi}{\to}S_{DVR}$.
\epr
\beR
As one sees from the end of this proof, we have a stronger property: $\overline{ann\big(T^1_{(\Si,G_{lr},A)}\big)\cdot I_{m-1}(A)}\sseteq\overline{I_m(A)}$.

One would like to strengthen it further, to achieve from $\phi(f)\in\phi\big(I_m(A)\big):\phi\big(I_{m-1}(A)\big)$ that $\phi(f)\in\phi\big(I_m(A):I_{m-1}(A)\big)$. Then the final conclusion would be: $ann\big(T^1_{(\Mat,G_{lr},A)}\big)\sseteq \overline{I_m(A):I_{m-1}(A)}$.
But in general: $\phi(I:J)\neq\phi(I):\phi(J)$.
\eeR

\bex\label{Ex.ann(T1).radical.ideal}
Let $R$ be a  Noetherian ring and $A\in Mat_{m\times n}(R)$.  Suppose that the ideal $I_m(A)$ is radical  and $height(I_{m-1}(A))>height(I_{m}(A))$.
Then $ann.coker(A)=(det(A))$, see Property \eqref{Thm.Background.ann.coker.in.terms.of.Fitting.ideals}. In addition $I_m(A)=\overline{I_m(A)}$.
 Finally, as $height(I_{m-1}(A))>height(I_{m}(A))$, we have
$\overline{I_m(A)}=\overline{I_m(A)}:\overline{I_{m-1}(A)}$.
 Altogether we get:  $ann(T^1_{(\Si,G_{lr},A)})=(det(A))$.
\eex
\bex\label{Ex.ann(T1).over.PID}
Suppose $R$ is a PID. Then all ideals are principal, hence integrally closed. Thus, for any $A\in\Mat$ we have
\beq
ann.coker(A)=I_m(A):I_{m-1}(A)=\overline{I_m(A)}:\overline{I_{m-1}(A)}.
 \eeq
 (For example, bring $A$ to the Smith normal form.) Therefore  $ann(T^1_{(\Si,G_{lr},A)})=I_m(A):I_{m-1}(A)=ann.coker(A)$.
\eex

\subsubsection{Applications to finite determinacy}

\

{\em Proof of Proposition \ref{Thm.Results.Fin.Det.for.GlGrGlr}.}

We use Proposition \ref{Thm.Intro.ord.of.det.via.Loewy.length}. First we remark that  the assumption \eqref{Eq.Intro.assumption.R.Noetherian.or.surjects} ensures the needed approximation property of $R$, see \S\ref{Sec.Background.Approximation.Properties}.
Indeed,  in all the cases the equations $UA=A+B$, $AV^{-1}=A+B$, $UAV^{-1}=A+B$ are linear in the entries of $U,V$.
 (In the later case present the equation in the form $UA=(A+B)V$.)

{\bf 1.} Note that $ann(T^1_{(\Si,G^{(m)}_r,A)})=ann\quots{R^{\oplus
m}}{\cm\cdot Im(A)}= ann\quots{R^{\oplus m}}{Im(A|_{\cm\cdot R^{\oplus
n}})}$.
 Thus Proposition \ref{Thm.Intro.ord.of.det.via.Loewy.length}
gives the stated bound on $ord^\Si_{G_r}(A)$.

{\bf i.}  If $dim(R)>(n-m+1)$ then for any $A\in Mat_{m\times n}(\cm)$ the ideal $ann.coker(A)$,
having height at most $(n-m+1)$,  see Property (\ref{Thm.Background.ann.coker.in.terms.of.Fitting.ideals}),
cannot contain any power of the maximal ideal $\cm$.
This implies the statement.

{\bf ii.} Vice versa,  the height of the ideal  $ann.coker(A)$ is generically the expected one, see \S\ref{Sec.Background.Fitting.Ideals.ann.coker}.
Thus for any $A$ and the generic $B\in Mat_{m\times n}(\cm^N)$, the ideal $ann.coker(A+B)$
is of height $min\Big((n-m+1),dim(R)\Big)$. Thus if $dim(R)\le (n-m+1)$ then $ann.coker(A+B)$
contains a (finite)  power of the maximal ideal.

{\bf 2.} Follows immediately from part 2 of Theorem \ref{Thm.Results.Annihilators.for.GlGrGlr}.

{\bf 3. i.} Note that $height\big(\overline{I_m(A)}:\overline{I_{m-1}(A)}\big)=height(\overline{I_m(A)})=height(I_m(A))=height(ann.coker(A))$.
Thus the statement follows immediately from  the first part and part 3 of Theorem \ref{Thm.Results.Annihilators.for.GlGrGlr}.

{\bf 3.ii.}
Note that $G_{lr}\supset G_r$, thus, by part 1:
\beq
ord^\Si_{G_{lr}}(A)\le ord^\Si_{G_r}(A)\le ll_R(ann.coker(A|_{\cm R^p}))-1.
\eeq

If $R$  is Noetherian then the lower bound on $ord^\Si_{G_{lr}}(A)$ is immediate.

To prove the lower bound in the non-Noetherian case we
use the surjectivity of the completion, $R\stackrel{\hat{\psi}}{\twoheadrightarrow}\hR$.
Then by Lemma
\ref{Thm.Background.Completion.Lemma} we have the
equality of Loewy lengths:
\beq
ll_R(ann(T^1_{(\Si,G_{lr},A)}))=ll_\hR\Big(\hat{\psi}(ann(T^1_{(\Si,G_{lr},A)}))\Big)=
ll_\hR\Big(ann\quots{\hat{\psi}(T_{(\Si,A)})}{\hat{\psi}(T_{(G_{lr}A,A)})}\Big)=
ll_\hR\Big(\overline{I_m(\hA)}:\overline{I_{m-1}(\hA)}\Big). \eeq
\epr

\bex
 Consider the trivial case: $A$ is a ``numerical" matrix (with entries in $\k$), $m\le n$.
Then  $ann.coker(A)$ is either $R$ or $\{0\}$. Hence
$A$ is finitely $G_{lr}$-determined
iff at least one of its maximal minors is a non-zero constant, i.e. $A$ is of the full rank.
(In other words, $A$ is invertible from the left, i.e. it has no left-kernel.) In this case, for $m\le n$,
$A$ is 0-determined \wrt $G_r$.
\eex
\bex
Let $R=\k[[x,y]]$ and $A=\bpm x^5&0&y^3\\0&y^4&x^3\epm$. Then $I_1(A)=(x^3,y^3)$, $I_2(A)=(y^7,x^5y^4,x^8)$
and  $ann.coker(A)=(y^7,x^5y^4,x^8)$. (To compute $ann.coker$ one can notice that $height(I_2(A))=2$, as expected,
 thus $ann.coker(A)=I_2(A)$, by Property \eqref{Thm.Background.ann.coker.in.terms.of.Fitting.ideals}.)
Therefore $\overline{I_1(A)}=\cm^3$ and  $\overline{I_2(A)}=\cm^8+(y^7)$.
Thus
\beq
(y^7,x^5y^4,x^8)\sseteq ann T^1_{(\Si,G_{lr},A)}\sseteq \overline{I_2(A)}:\overline{I_1(A)}=\cm^5.
\eeq
Note that $ann.coker(A|_{\cm\cdot R^{\oplus 3}})=(y^7,x^5y^4,x^8)$.
Finally, $ll_R(\cm^5)=5$ and $ll_R(y^7,x^5y^4,x^8)=11$.
We get: $4\le ord^\Si_{G_{lr}}(A)\le 10$.
\eex
\bex\label{Thm.Finite.Determinacy.Glr.Low.Dimensional.Cases} Suppose  $R$ satisfies the assumption \eqref{Eq.Intro.assumption.R.Noetherian.or.surjects} and $dim(R)>0$.  Let $\Si=\Mat$.

1. For $m=n$,  $A\in Mat_{m\times m}(\cm)$ is finitely-$G_{lr}$-determined iff $dim(R)=1$, $det(A)\in R$ is not a
zero divisor and $det(A)\not\in\cm^\infty$. This generalizes Part 1 of \cite[Theorem 1.1]{Bruce-Tari04}.

Proof: for square matrices $height\Big(ann.coker(A)\Big)=height(det(A))\le 1$. Thus it can contain a power of the maximal ideal only when $dim(R)\le1$. Now invoke Proposition \ref{Thm.Results.Fin.Det.for.GlGrGlr}.

2. Suppose $dim(R)=2$ and $m<n$. Suppose $A$ has at least two $m\times m$ blocks whose determinants, $\De_1,\De_2$, are neither zero divisors,
 nor belong to $\cm^\infty$ and moreover are relatively prime,
 i.e. if $\De_i=a_i h\in R$ then $h\in R$ is invertible.  Then A is finitely $G_r$-determined.

Proof: Let $\De_1,\De_2$ be two such determinants then the height of the ideal $(\De_1)+(\De_2)$ is two.
 Thus $height(ann.coker(A))=2$ and $ann.coker(A)$ contains a power of $\cm$.
\eex
\bex
 Suppose  $R$ satisfies the assumption \eqref{Eq.Intro.assumption.R.Noetherian.or.surjects} and $A\in\Mat$.
Then part 3 of Theorem \ref{Thm.Results.Annihilators.for.GlGrGlr}
 and Theorem \ref{Thm.Background.Linearization} imply:
\beq
\text{if $\tA\equiv A\ mod(\cm\cdot ann.coker(A))$ then $\tA\stackrel{G_{lr}}{\sim}A$}.
\eeq
   This both strengthens (quantifies) and generalizes  \cite[Theorem 5.2]{Cutkosky-Srinivasan}.
\eex

\subsubsection{Applications to the infinite determinacy}
Part 1 of Proposition \ref{Thm.Results.Fin.Det.for.GlGrGlr} implies the following criterion.
\bcor\label{Thm.Proofs.Infinite.Determinacy.Gr}
Let $A\in Mat_{m\times n}(C^\infty(\R^p,0))$. If there exists an element $f\in I_m(A)$ satisfying {\L}ojasiewicz inequality, Equation \eqref{Eq.Lojasiewicz.Ineq}, then $A$ is infinitely-$G_r$-determined.
\ecor
Indeed, by Proposition \ref{Thm.Results.Fin.Det.for.GlGrGlr} it is enough to check that the ideal $I_m(A)$ contains $\cm^\infty$, in other words: any function flat at the origin is divisible by some element of $I_m(A)$. And this is ensured by {\L}ojasiewicz inequality.

\subsection{(Anti-)symmetric matrices and the congruence}\label{Sec.Proofs.Congruence}

\subsubsection{Computation of $ann(T^1_{(\Si,G_{congr},A)})$} \
{\em Proof of Theorem \ref{Thm.Results.Annihilator.T1.for.Gcongr}.}

{\bf 1.}  First we assume that $R$ is Noetherian.
To prove that $ann(T^1_{(\Si,G_{congr},A)})=\{0\}$ it is enough to show that this module is not a torsion module and is supported at the generic point
 of any component of $Spec(R)$. Thus we assume that $Spec(R)$ is irreducible and localize at the generic point of $Spec(R)$. (Note that $dim(R)>0$.) At this
 point we compare the modules, or rather vector spaces, $(T_{(GA,A)})_{(0)}$ and $(T_{(\Si,A)})_{(0)}$.

Note that $R_{(0)}$ is a field, thus we can bring the symmetric part of $A$ to the canonical form of
 Proposition \ref{Thm.Background.Canonical.Form.Matrix.over.DVR}, i.e.
$A\stackrel{G_{congr}}{\sim} Diag+A_-$, where $Diag$ is diagonal while $A_-\in Mat^{skew-sym}_{m\times m}(R)$.
We can assume that $A+A^t$ is generically non-degenerate on $Spec(R)$, therefore $Diag$ is invertible over $R_{(0)}$.
To understand  $(T_{(GA,A)})_{(0)}$, which is spanned by $UA+AU^t$, we re-scale the matrices, $\tU:= U\cdot Diag$, and expand them into
 the symmetric and skew-symmetric  parts, $\tU=\tU_++\tU_-$.
Then
 $(T_{(GA,A)})_{(0)}=Span_R(2\tU_+,\tU A_-+A_-\tU^t)$. To prove that
\beq
(T_{(GA,A)})_{(0)}\subsetneq(T_{(\Si,A)})_{(0)}=Mat^{sym}_{m\times m}(R_{(0)})\oplus Mat^{skew-sym}_{m\times m}(R_{(0)})
\eeq
 we consider the map
\beq
 Mat_{m\times m}(R_{(0)})\stackrel{\phi}{\to}Mat^{sym}_{m\times m}(R_{(0)})\oplus Mat^{skew-sym}_{m\times m}(R_{(0)}),\quad
 \tU\to(2\tU_+,\tU A_-+A_-\tU^t)
\eeq
Note that the domain and the target are $R_{(0)}$-vector spaces of the same (finite) dimension. Thus it is enough to demonstrate
a non-trivial kernel of $\phi$. Indeed, the space $ker(\phi)=\{\tU_+=\zero,\ \tU_- A_-=A_-\tU_-\}$ is non-zero. It contains all the
 matrices $\tU_-$ that commute with $A_-$, e.g.
 $ker(\phi)\supseteq Span_R(A_-,A^3_-,A^5_-,\dots)$.

 \

Now, if $R$  is non-Noetherian, but satisfies condition \eqref{Eq.Intro.assumption.R.Noetherian.or.surjects}, we use Lemma \ref{Thm.Background.Completion.Lemma} to get: $\hat\psi(ann(T^1))=\{0\}$. Therefore $ann(T^1)\sseteq\cm^\infty$.

{\bf 2.} To prove that $ann(T^1_{(\Si,G_{congr},A)}) \supseteq ann.coker(A)$, take any $B\in Mat^{sym}_{m\times m}(ann.coker(A))$. By the definition of the annihilator-of-cokernel there exists $U$ satisfying: $AU=\frac{B}{2}$. As $A$, $B$ are symmetric: $U^tA=\frac{B}{2}$. Therefore $UA+AU^t=B$.
 But $UA+AU^t\in T_{(G_{congr}A,A)}$, see Example \ref{Ex.Tangent.Spaces.to.the.Actions}.
 Thus $ann.coker(A)\cdot T_{(\Si,A)}=Mat^{sym}_{m\times m}(ann.coker(A))\sseteq T_{(G_{congr}A,A)}$.

To prove that $ann(T^1_{(\Si,G_{congr},A)})\sseteq\overline{I_m(A)}:\overline{I_{m-1}(A)}$, let $f\in ann(T^1_{(\Si,G_{congr},A)})$, i.e. $f\cdot Mat^{sym}_{m\times m}(R)\sseteq Span_R(UA+AU^t|\ U\in Mat_{m\times m}(R)$. Then, for any $R\stackrel{\phi}{\to}S_{DVR}$:
\beq
\phi(f)\cdot Mat^{sym}_{m\times m}(S)\sseteq Span_R(U\phi(A)+\phi(A)U^t|\ U\in Mat_{m\times m}(S)).
\eeq
 Now, bring $\phi(A)$ to the canonical form (using Proposition \ref{Thm.Background.Canonical.Form.Matrix.over.DVR}),  $\phi(A)=V\cdot D_{\phi(A)}V^t$, where $D_{\phi(A)}=\oplus\la_i\one$, with $(\la_1)\supseteq(\la_2)\supseteq\cdots\supseteq(\la_m)$.
 Therefore:
\beq
\phi(f)\cdot Mat^{sym}_{m\times m}(S)\sseteq Span_R(UD_{\phi(A)}+D_{\phi(A)}U^t|\ U\in Mat_{m\times m}(S)).
 \eeq
 But the $(m,m)$-entry of any matrix of the form $(UD_{\phi(A)}+D_{\phi(A)}U^t)$  lies in the ideal
 \beq
 (\la_m)=ann.coker(\phi(A))=I_m(\phi(A)):I_{m-1}(\phi(A))=\phi(I_m(A)):\phi(I_{m-1}(A)).
 \eeq
  Therefore
$\phi(f\cdot I_{m-1}(A))\in \phi(I_m(A))$.
As this holds for any $R\stackrel{\phi}{\to}S_{DVR}$ we get:
$\overline{f\cdot I_{m-1}(A)}\sseteq\overline{I_m(A)}$. Finally, Equation (\ref{Eq.Background.integral.closure.product.of.ideals}) implies
 $f\in \overline{I_{m}(A)}:\overline{I_{m-1}(A)}$.

{\bf 3.} First we prove the part $(\cdots)\sseteq ann(T^1_{(\Si,G_{congr},A)})$.

\quad 3.i. For $m$-even the proof of $ann(T^1_{(Mat_{m\times m}^{skew-sym}(R),G_{congr},A)})\supseteq ann.coker(A)$ goes as in part 2.

\quad 3.ii.   Suppose $m$ is odd, then by Proposition \ref{Thm.Background.AntiSymmetric.Matrices.Pfaffians}:
 $T_{(G_{congr}A,A)}\supseteq Mat^{skew-sym}_{m\times m}(Pf_{m-1}(A))$.
 Thus
 $ann(T^1_{(\Si,G_{congr},A)})\supseteq (Pf_{m-1}(A))$.

\

To establish the bound   $ann(T^1_{(\Si,G_{congr},A)})\sseteq(\cdots)$, let $f\in ann(T^1_{\Si,G_{congr},A})$, i.e. $f\cdot  Mat^{skew-sym}_{m\times m}(R)\sseteq Span_R(UA+AU^t|\ U\in Mat_{m\times m}(R)$. Then, for any $R\stackrel{\phi}{\to}S_{DVR}$:
\beq
 \phi(f)\cdot Mat^{skew-sym}_{m\times m}(S)\sseteq Span_R(U\phi(A)+\phi(A)U^t|\ U\in Mat_{m\times m}(S)).
\eeq
Bring $\phi(A)$ to the canonical form (using Proposition \ref{Thm.Background.Canonical.Form.Matrix.over.DVR}), $\phi(A)=V\cdot E_{\phi(A)}V^t$, where $E_{\phi(A)}=\oplus_i \la_i E_i$
   and $(\la_1)\supseteq(\la_2)\cdots\supseteq(\la_{\lfloor\frac{m}{2}\rfloor})$.
   Therefore:
\beq
\phi(f)\cdot Mat^{skew-sym}_{m\times m}(S)\sseteq Span_R(UE_{\phi(A)}+E_{\phi(A)}U^t|\ U\in Mat_{m\times m}(S)).
\eeq
The $(m-1,m)$ and $(m,m-1)$ entries of $(UE_{\phi(A)}+E_{\phi(A)}U^t)$ belong to the ideal
   $(\la_{\lfloor\frac{m}{2}\rfloor})$.
Therefore,
\beq\label{Eq.inside.proof}
\phi(f)\in (\la_{\lfloor\frac{m}{2}\rfloor})=I_{m-1}(\phi(A)):I_{m-2}(\phi(A))=\phi(I_{m-1}(A)):\phi(I_{m-2}(A)), \eeq
i.e. $\phi(f\cdot I_{m-2}(A))\sseteq\phi(I_{m-1}(A))$.

As this holds for any $R\stackrel{\phi}{\to}S_{DVR}$ we get:
$\overline{f\cdot I_{m-2}(A)}\sseteq\overline{I_{m-1}(A)}$. Finally, we get
 $f\in \overline{I_{m-1}(A)}:\overline{I_{m-2}(A)}$, see Equation (\ref{Eq.Background.integral.closure.product.of.ideals}).
Therefore $ann(T^1_{(\Si,G_{congr},A)})\sseteq \overline{I_{m-1}(A)}:\overline{I_{m-2}(A)}$.

\

Finally, if $m$ is even, then $\phi(I_{m}(A)):\phi(I_{m-1}(A))=\phi(I_{m-1}(A)):\phi(I_{m-2}(A))$. Thus
in Equation (\ref{Eq.inside.proof}) we have: $\phi(f)\in \phi(I_{m}(A)):\phi(I_{m-1}(A))$ and therefore $ann(T^1_{(Mat_{m\times m}^{skew-sym}(R),G_{congr},A)})\sseteq \overline{I_{m}(A)}:\overline{I_{m-1}(A)}$.
\epr
\bex
(Continuing examples \ref{Ex.ann(T1).radical.ideal} and \ref{Ex.ann(T1).over.PID}.) In some cases the bounds on $ann(T^1_{(\Si,G_{lr},A)})$ turn into equality.
Let $R$ be Noetherian ring and $A$ either symmetric or skew-symmetric of even size.  Suppose
\bee[i.]
\item either $R$ is PID
\item or $I_m(A)$ is radical  and $height(I_{m-1}(A))>height(I_{m}(A))>0$.
\eee
 Then $ann(T^1_{(\Si,G_{lr},A)})=ann.coker(A)$. Indeed, $\overline{I_m(A)}=I_m(A)$ (Property \ref{Thm.Background.Integral.Closure.Properties}) and
  $ann.coker(A)=I_m(A):I_{m-1}(A)$ (Property \ref{Thm.Background.Ann.Coker.Properties}).
\eex

\bex
In the low dimensional cases $T^1_{(\Si,G_{congr},A)}$ is easily written down explicitly.
\\1. In the trivial case $A=\bpm 0&a\\-a&0\epm$ we get $ann.coker(A)=(a)\sset R$ and thus
$T^1_{(\Si,G_{congr},A)}\approx(\quots{R}{(a)})^{\oplus 2}$.
\\2. For $3\times 3$ matrices take $A=\bpm 0&a&b\\-a&0&c\\-b&-c&0\epm$, matrix of indeterminates, we get
\[
(Pf_2(A))=(a)+(b)+(c)=\overline{I_2(A)}:\overline{I_1(A)}\sset R.
\]
Therefore $T^1_{(\Si,G_{congr},A)}\approx\Big(\quots{R}{(a)+(b)+(c)}\Big)^{\oplus 3}$.
\eex

\subsubsection{Finite determinacy for congruence}

\

{\em Proof} of Proposition \ref{Thm.Results.Finite.Determinacy.G_congr}.
The statements 1, 1', 1'' follow from Theorem \ref{Thm.Results.Annihilator.T1.for.Gcongr} because $height(ann.coker(A))\le1$ for a square matrix,
  while for a skew-symmetric matrix of odd-size $height(Pf_{m-1}(A))\le 3$, see \S\ref{Sec.Background.Fitting.Ideals.ann.coker}.

The statements 2.i. 2.ii. follow from Theorem \ref{Thm.Results.Annihilator.T1.for.Gcongr} and Proposition \ref{Thm.Intro.ord.of.det.via.Loewy.length}.

\bex
The natural question is the order of determinacy for a generic (skew-)symmetric matrix.
 By Proposition  \ref{Thm.Results.Finite.Determinacy.G_congr} we must assume $dim(R)=1$.
  For simplicity we assume that $R$ is a Henselian DVR.
\bee[i.]
\item Take the generic matrix $A\in Mat^{sym}_{m\times m}(\cm^k)$. Use Proposition \ref{Thm.Background.Canonical.Form.Matrix.over.DVR} to
bring  $A$ to the canonical form. (Note that the order of determinacy and $ann(T^1_{(\Si,G_{congr},A)})$ are invariant under the $G_{congr}$-action.)
 Thus $A=\oplusl_j \la_j\one$ and by genericity the orders of all $\la_j$ are $k$.  Thus  the order of
 determinacy is $k$. We get the rigidity in the following sense. Denote $\Si^{(k+1)}=\{A\}+Mat^{sym}_{m\times m}(\cm^{k+1})$,
  then $T^1_{(\Si^{(k+1)},G_{congr},A)}=\{0\}$.

\item Similarly for skew-symmetric matrices of even size. Let $A\in Mat^{skew-sym}_{m\times m}(\cm^k)$ be generic, then
  $A\stackrel{G_{congr}}{\sim}\oplusl_j \la_j E_j$, where $ord(\la_j)=k$. Thus the order of determinacy is $k$. If one takes
$\Si^{(k+1)}=\{A\}+Mat^{skew-sym}_{m\times m}(\cm^{k+1})$ then $T^1_{(\Si^{(k+1)},G_{congr},A)}=\{0\}$.
\eee\eex

\

 Proposition \ref{Thm.Results.Finite.Determinacy.G_congr} implies the following criterion for the infinite determinacy of (skew-)symmetric forms
 valued in the germs of smooth functions.
\bcor\label{Thm.Proofs.Infinite.Determinacy.G_congr}
1. Let either $A\in Mat^{sym}_{m\times m}(C^\infty(\R^p,0))$ or $A\in Mat^{skew-sym}_{m\times m}(C^\infty(\R^p,0))$, for $m$-even. If $det(A)$ satisfies the {\L}ojasiewicz inequality,
 Equation \eqref{Eq.Lojasiewicz.Ineq},
 then $A$ is infinitely-$G_{congr}$-determined.
\\2. Let $A\in Mat^{skew-sym}_{m\times m}(C^\infty(\R^p,0))$, for $m$-odd. If there exists an element $f\in I_{m-1}(A)$ that satisfies the {\L}ojasiewicz inequality, then $A$ is infinitely-$G_{congr}$-determined.
\ecor
The argument is the same as for Corollary \ref{Thm.Proofs.Infinite.Determinacy.Gr}.

\subsection{The upper-block-triangular matrices and the action of $G^{up}_r$, $G^{up}_{lr}$}\label{Sec.Proofs.Upper.Triangular}
First we give two examples showing that the block-triangular matrices appear naturally, then we prove Theorem \ref{Thm.Results.Annihilator.T1.Upper.Triang}.
\subsubsection{Morphisms of filtered free modules and their deformations}
Let $M$ be a free module and denote by $M_\bullet$  a filtration by free submodules which are direct summands, $\{0\}=M_0\ssetneq M_1\ssetneq \cdots\ssetneq  M_k=M$.
 Given two such filtered modules, $M_\bullet$, $N_\bullet$, consider their filtered homomorphisms,
 $Hom_R(M_\bullet,N_\bullet)=\{\phi\in Hom_R(M,N)|\ \phi(M_i)\sseteq N_i\}$. Fix a filtered basis of $M$, i.e. a minimal set of generators
  of $M$ compatible with the filtration:
 \beq
 (e_1,\dots,e_{i_1},e_{i_1+1},\dots,e_{i_2},\dots,e_{i_{k-1}+1},\dots,e_{i_k}), \quad \Big\{M_j=Span_R(e_1,\dots,e_{i_j})\Big\}_{j=1,\dots,k}.
\eeq
Similarly fix an $N_\bullet$-basis. For these chosen bases any homomorphism $\phi\in Hom_R(M_\bullet,N_\bullet)$ is presented by an upper-block-triangular matrix,
\beq
[\phi]=A=\bpm A_{11}&A_{12}&\dots\dots&\dots&A_{1k}\\\zero&A_{22}&\dots&\dots&A_{2k}\\\dots&\dots&\dots&\dots&\dots\\\zero&\zero&\dots&\zero&A_{kk}\epm.
\eeq
The change of bases in $M_\bullet$, $N_\bullet$ act on $A$ by elements of $G^{up}_{lr}=G^{up}_l\times G^{up}_r
=GL(M_\bullet,R)\times GL(N_\bullet,R)$.
Thus for the deformations of $\phi$ inside $\Si=Hom_R(M_\bullet,N_\bullet)$ we have:
\beq
T^1_{(\Si,G^{up}_{lr},A)}= T^1_{\phi\in Hom_R(M_\bullet,N_\bullet)}=\quotient{T_{(Mat^{up},A)}}{T_{(G^{up}_{lr}A,A)}}.
\eeq

\subsubsection{Deformations of filtered submodules}
Fix a filtered module $\{0\}=M_0\ssetneq M_1\ssetneq \cdots\ssetneq  M_k=M\sseteq R^{\oplus m}$, not necessarily free.
 Suppose a sequence of ``containers" is prescribed, i.e. $M_i\sseteq V_i$, where $V_i\sseteq R^{\oplus m}$ is a free direct summand.
   We assume $V_i\subsetneq V_{i+1}$. Choose a set
   of generators of $R^{\oplus m}$ compatibly with the filtration $V_i\ssetneq V_{i+1}\ssetneq\dots R^{\oplus m}$.
Then the generating matrix $A_M$ of $M$ is upper-block-triangular:
\beq
A_M=\bpm A_{11}&A_{12}&\dots\dots&\dots&A_{1k}\\\zero&A_{22}&\dots&\dots&A_{2k}\\\dots&\dots&\dots&\dots&\dots\\\zero&\zero&\dots&\zero&A_{kk}\epm.
\eeq
A change of basis of $M_\bullet$ corresponds to the $G^{up}_r$-action. Therefore for the deformations of $M_\bullet$, inside the ambient filtration $V_\bullet$,
 we have
\beq
T^1_{(\Si,G_r,A)}=T^1_{M_\bullet, V_\bullet}=\quotient{T_{(Mat^{up},A)}}{T_{(G^{up}_{r}A,A)}}.
\eeq

\subsubsection{}
{\em Proof of Theorem \ref{Thm.Results.Annihilator.T1.Upper.Triang}.}

{\bf 1.} Denote $\Si=Mat^{up}_{\{m_i\}\times\{n_j\}}(R)$, thus $T_{(\Si,A)}\approx Mat^{up}_{\{m_i\}\times\{n_j\}}(R)$. We write down the tangent space to the group orbit:
\beq
T_{(G^{up}_rA,A)}=Span_R\Bigg(\bpm A_{11}U_{11}&A_{11}U_{12}+A_{12}U_{22}&\dots&\dots \suml^k_{i=1}A_{1i}U_{ik}\\
\zero&A_{22}U_{22}&\dots&\suml^k_{i=2}A_{2i}U_{ik}\\\dots&\dots&\dots&\dots&
\\
\zero&\dots&\zero&A_{kk}U_{kk}
\epm,\ \Big\{U_{ij}\in Mat_{m_i\times n_j}(R)\Big\}_{i\le j}\Bigg)
\eeq
Note that both modules $T_{(\Si,A)}$ and $T_{(G^{up}_rA,A)}$ decompose into the direct sums, according to the blocks of columns. Explicitly, $T_{(\Si,A)}=\oplusl^k_{l=1}\La_l$, where $\La_l$ consists of all matrices $B$ with the block-structure $B=\{B_{ij}\}$ satisfying: $B_{ij}=\zero$ if $j\neq l$ or if $i>j$.

Therefore, for a fixed $A$ we have:
$T^1_{(\Si,G^{up}_r,A)}=\oplusl^k_{q=1}(T^1_{(\Si,G^{up}_r,A)})_q$, where
\begin{multline}
(T^1_{(\Si,G^{up}_r,A)})_q=\\
\frac{Mat_{(\suml^q_{j=1}m_j)\times n_q}(R)}{Span_R\Bigg(\bpm A_{11}&\dots&\dots&A_{1q}\\\zero&\dots&\dots&\\
\zero&\dots&\zero& A_{qq}\epm\bpm U_{1q}\\\dots\\U_{qq}\epm,\ U_{jq}
\in Mat_{m_j\times n_q}(R)\Bigg)}=
coker\bpm A_{11}&\dots&\dots&A_{1q}\\\zero&\dots&\dots&\\\zero&\dots&\zero&A_{qq}\epm=ann.coker(\cornA{q}).
\end{multline}
 Accordingly: $ann(T^1_{(\Si,G^{up}_r,A)})=\capl^k_{q=1}ann(T^1_{(\Si,G^{up}_r,A)})_q=\capl^k_{q=1}ann.coker(\cornA{q})$.
This proves the first statement.

\

{\bf 2.}
The bounds $ \prodl^k_{q=1} ann.coker(A_{qq})\sseteq \capl^k_{q=1}ann.coker(\cornA{q})\sseteq \capl^k_{q=1} ann.coker(A_{qq})$ are proved by induction as follows.
The case $k=1$ is trivial. Suppose the bounds hold for $(k-1)$. To establish them for $k$ it is enough to prove:
\bee[i.]
\item
 $ann.coker(A)\sseteq ann.coker(A_{kk})$
\item
 $ann.coker(A)\supseteq \prodl^k_{q=1}ann.coker(A_{qq})$.
\eee
In the proof we decompose $R^{\sum m_q}=\oplusl_q R^{ m_q}$. To keep track of the block structure we denote the image of $R^{ m_q}$
 inside this sum by $R^{ m_q}\oplus\{0\}$.

\quad {\em Proof of i.}
 If $f\in ann.coker(A)$ then $f\cdot(\{0\}\oplus R^{ m_k})\sseteq Im(A)$ and thus $f\cdot R^{ m_k}\sseteq Im(A_{kk})$, i.e. $f\in ann.coker(A_{kk})$.

\quad {\em Proof of ii.} It is enough to check: for any collection $\{f_q\in ann.coker(A_{qq})\}_{q=1\dots k}$ holds: $\prodl^k_{q=1} f_q\in ann.coker(A)$.

By the inductive assumption we have:
$(\prodl^{k-1}_{q=1} f_q)\cdot(\oplusl^{k-1}_{q=1}  R^{m_q}\oplus\{0\})\sseteq Im(A)$. It remains to check:
 $(\prodl^{k}_{q=1} f_q)\cdot(\{0\}\oplus  R^{m_k})\sseteq Im(A)$. Indeed,  any $0\oplus v\in (\prodl^{k}_{q=1} f_q)\cdot(\{0\}\oplus  R^{m_k})$ is presentable
  as $v=A_{kk}u$ for some $u\in (\prodl^{k-1}_{q=1} f_q)R^{n_k}$.
  Thus $v-A(0\oplus u)\in (\prodl^{k-1}_{q=1} f_q)\cdot(\oplusl^{k-1}_q  R^{m_q}\oplus\{0\})\sseteq Im(A)$.

\

{\bf 3.} The lower bound on $ann(T^1_{(\Si,G^{up}_{lr},A)})$ follows by parts 1,2 because $G^{up}_{lr}\supset G^{up}_r$.
 For the upper bound we write down the tangent space
\beq\small
T_{(G^{up}_{lr}A,A)}=Span_R\Bigg(\bpm A_{11}U_{11}&\dots& \suml^k_{i=1}A_{1i}U_{ik}\\
\zero&\dots&\suml^k_{i=2}A_{2i}U_{ik}\\\dots&\dots&\dots
\\
\zero&\zero&A_{kk}U_{kk}\epm,
\bpm V_{11}A_{11}&\dots& \suml^k_{i=1}V_{1i}A_{ik}\\
\zero&\dots&\suml^k_{i=2}V_{2i}A_{ik}\\\dots&\dots&\dots
\\
\zero&\zero&V_{kk}A_{kk}\epm\Bigg)
\eeq
Thus $f\cdot T_{(\Si,A)}\sseteq T_{(G^{up}_{lr}A,A)}$ implies $Mat_{m_i\times n_i}((f))\sseteq Span_R(V_{ii}A_{ii}+A_{ii}U_{ii})$ for any $i$. But the later means:
\beq
f\in \capl_i ann(T^1_{(\Si_i,G_{lr},A_{ii})})\sseteq\capl_i \Big(\overline{I_{m_i}(A_{ii})}:\overline{I_{m_i-1}(A_{ii})}\Big).\text{\epr}
\eeq

\bex
\bee[i.]
\item Suppose the blocks $\{A_{qq}\}$ are mutually generic so that the ideals $ann.coker(A_{qq})$ are relatively prime in the sense:
 $\prodl^{k}_{q=1}ann.coker(A_{qq})=\capl^{k}_{q=1}ann.coker(A_{qq})$.
 Then  we get the equality, $ann(T^1_{(\Si,G^{up}_r,A)})=\prodl^{k}_{q=1}ann.coker(A_{qq})$.
\item Suppose $A_{ij}=\zero$  for $1\le i<j\le k$. Then we have the obvious $ann(T^1_{(\Si,G^{up}_r,A)})=\capl^{k}_{i=1}ann.coker(A_{ii})$.
\item Consider the case of one Jordan cell,
\[
A=\bpm a&1&0&\dots&\dots\\0&a&1&0&\dots\\\dots&\dots&\dots&\dots&\dots\\0&\dots&\dots&0&a\epm.
\]
Then, by direct check, $ann.coker(A)=(a^n)$. This realizes the lower bound.
\item Suppose $A$ is a square matrix, hence $\{m_q=n_q\}$, and $det(A)$ is not a zero divisor in $R$. Suppose, moreover, $height(I_{m_q-1}(A_{qq}))>1$
 for all the blocks. Then we get $(det(A))\sseteq ann(T^1_{(\Si,G^{up}_r,A)})\sseteq \capl_q(det(A_{qq}))$.
\eee
\eex

\subsubsection{Applications to finite determinacy}
 Theorem \ref{Thm.Results.Annihilator.T1.Upper.Triang} implies the immediate
\bcor
\bee
\item If $dim(R)>(n_i-m_i+1)$  for some $i$, then no matrix $A\in Mat^{up}_{\{m_i\}\times\{n_j\}}(\cm)$ is finitely-$G^{up}_{lr}$-determined.
\item Suppose $R$ satisfies condition  \eqref{Eq.Intro.assumption.R.Noetherian.or.surjects} and moreover $dim(R)\ge(n_i-m_i+1)$ for all $i$.
 Then the generic finite-$G^{up}_{r}$-determinacy holds and
\begin{multline}
ll_R\Big(\capl_i ann.coker(A_{ii})\Big)-1\le ll_R\Big(\capl_q ann.coker(\cornA{q})\Big)\le ord^{\Si}_{G^{up}_r}
\le\\\le ll_R\Big(\capl_q ann.coker(\cornA{q}|_{\cm\cdot R^{\suml^q_{j=1} n_j}})\Big)
\le
ll_R\Big(\prod_i ann.coker(A_{ii}|_{\cm\cdot R^{\oplus m}})\Big)-1.
\end{multline}
\eee
\ecor

\subsection{The case of conjugation}
{\em Proof of Theorem \ref{Thm.Results.T1.Conjugation}.}  Denote $\Si=Mat_{m\times m}(R)$.

1. First we assume that $R$ is Noetherian.
As in the case of congruence, to prove that $ann(T^1_{(\Si,G_{conj},A)})=\{0\}$ it is enough to prove that $T^1_{(\Si,G_{conj},A)}$
is not a torsion on any irreducible component of $Spec(R)$. Thus we restrict to an irreducible component of $Spec(R)$ and check $T^1_{(\Si,G_{conj},A)}$
 at the generic point of $Spec(R)$, algebraically  we localize at the generic ideal $\{0\}$.
 Therefore we compare the $R_{(0)}$-vector spaces:
\beq
Mat_{m\times m}(R_{(0)})\quad vs \quad (T_{(G_{conj}A,A)})_{(0)}=Span_{R_{(0)}}\Big(UA-AU|\ U\in Mat_{m\times m}(R_{(0)})\Big).
\eeq
Now the inequality $(T_{(G_{conj}A,A)})_{(0)}\ssetneq (T_{(\Si,A)})_{(0)}$ is immediate as $(T_{(G_{conj}A,A)})_{(0)}$ lies inside the hyperplane
 $\{B|\ trace(B)=0\}\sset (T_{(\Si,A)})_{(0)}$.
 Therefore $\Big(T^1_{(\Si,G_{conj},A)}\Big)_{(0)}\neq\{0\}$.

If $R$ is not Noetherian but satisfies condition \eqref{Eq.Intro.assumption.R.Noetherian.or.surjects},
then we
take the completion  $R\stackrel{\hat\psi}{\to}\hR$. We have $ann(\hat\psi(T^1_{(\Si,G_{conj},A)}))=\{0\}$,
thus by Lemma \ref{Thm.Background.Completion.Lemma.surjective.case}, $\hat\psi(ann(T^1_{(\Si,G_{conj},A)}))=\{0\}$. Therefore
 (as in the case of congruence) we get $ann(T^1_{(\Si,G_{conj},A)}))\sseteq\cm^\infty$.

2. Follows immediately from part one and Proposition \ref{Thm.Intro.ord.of.det.via.Loewy.length}.
\epr

\subsection{Finite determinacy of chains of morphisms}\label{Sec.Proofs.Complexes}
Given a finite collection of free modules, $\{F_i\}_{i\in I}$, of ranks $\{m_i\}_{i\in I}$, and their morphisms, $\phi_i\in Hom_R(F_i,F_{i-1})$, we consider the corresponding chain $(F_\bullet,\phi_\bullet)$, not necessarily a complex. We assume $\{m_{i+1}\ge m_i\}$ and omit the terms $F_i$ with $m_i=0$.
The chains $(F_\bullet,\phi_\bullet)$, $(F_\bullet,\psi_\bullet)$ are called isomorphic if there exist automorphisms $U_i\circlearrowright F_i$ satisfying: $\psi_i\circ U_i=U_{i-1}\circ \phi_i$.
\bed
1. Denote by $\{J_{(F_i,\phi_i)}\}_{i\in I}$ the collection of the largest possible ideals satisfying:
\[\text{
if $\Big\{\psi_i\equiv\phi_i\ mod (J_{(F_i,\phi_i)})\Big\}_{i\in I}$,  then the chains $(F_\bullet,\phi_\bullet)$ and $(F_\bullet,\psi_\bullet)$ are isomorphic.}
\]

2. Denote by $J_{(F_\bullet,\phi_\bullet)}$ the largest possible ideal satisfying:
\[
\text{if $\Big\{\psi_i\equiv\phi_i\ mod (J_{(F_\bullet,\phi_\bullet)})\Big\}_{i\in I}$,  then the chains $(F_\bullet,\phi_\bullet)$ and $(F_\bullet,\psi_\bullet)$ are isomorphic.}
\]
\eed
\bel\label{Thm.Proofs.Lemma.Chains.Morphisms}
In both cases the largest object exists and  $J_{(F_\bullet,\phi_\bullet)}=\capl_i J_{(F_i,\phi_i)}$.
\eel
\bpr
(We prove the existence in the first case, the second case is similar.)

Call a collection of ideals $\{J_i\}_{i\in I}$ admissible if the conditions $\{\psi_i\equiv \phi_i\ mod(J_i)\}_{i\in I}$ imply the isomorphism $(F_\bullet,\psi_\bullet)\approx (F_\bullet,\phi_\bullet)$.
The collection of zero ideals, $\{J_i=\{0\}\}_{i\in I}$, is obviously admissible.
Thus it suffices to prove: if $\{J_i\}_{i\in I}$, $\{\tJ_i\}_{i\in I}$ are two admissible collections then the collection $\{J_i+\tJ_i\}_{i\in I}$ is admissible.

Suppose $\{\psi_i=\phi_i+\De_i+\tilde{\De}_i\}_{i\in I}$, where the entries of $\De_i$ lie in $J_i$, while the entries of $\tilde{\De}_i$ lie in $\tJ_i$. Then there exists a collection of automorphisms, $\{U_i\in Aut_R(F_i)\}_{i\in I}$ satisfying: $\{U_i\psi_iU^{-1}_i=\phi_i+U_i\tilde{\De}_iU^{-1}_i\}_{i\in I}$. Now, the entries of $U_i\tilde{\De}_iU^{-1}_i$ still lie in the ideal $\tJ_i$. Therefore there exists a collection of automorphisms, $\{V_i\in Aut_R(F_i)\}_{i\in I}$ satisfying: $\{V_iU_i\psi_iU^{-1}_iV^{-1}_i=\phi_i\}_{i\in I}$. The composition $\{V_iU_i\}_{i\in I}$ provides the needed automorphism of the chains.

\

Finally, both inclusions $J_{(F_\bullet,\phi_\bullet)}\sseteq\capl_i J_{(F_i,\phi_i)}$ and  $J_{(F_\bullet,\phi_\bullet)}\supseteq\capl_i J_{(F_i,\phi_i)}$ are obvious.
\epr

\bprop
Suppose $R$ is Noetherian.
\bee
\item $ann.coker(\phi_i|_{\cm\cdot F_i})\sseteq J_{(F_i,\phi_i)}\sseteq \overline{I_{m_i}(\phi_i)}:\overline{I_{m_i-1}(\phi_i)}$.
\item  $\capl_i ann.coker(\phi_i|_{\cm\cdot F_i})\sseteq J_{(F_\bullet,\phi_\bullet)}\sseteq \capl_i\Big(\overline{I_{m_i}(\phi_i)}:\overline{I_{m_i-1}(\phi_i)}\Big)$.
\eee
\eprop
This strengthens (quantifies) and generalizes \cite[Theorem 5.8]{Cutkosky-Srinivasan}.
\bpr
{\bf 1.} The isomorphism of the chains $(F_\bullet,\phi_\bullet)$, $(F_\bullet,\psi_\bullet)$ implies in particular the left-right equivalence of the representing matrices: $\phi_i\stackrel{G_{lr}}{\sim}\psi_i$ for each $i$. Thus the embedding $J_{(F_i,\phi_i)}\sseteq \overline{I_{m_i}(\phi_i)}:\overline{I_{m_i-1}(\phi_i)}$ follows immediately from part 3 of Theorem \ref{Thm.Results.Annihilators.for.GlGrGlr} and Theorem \ref{Thm.Background.Linearization}.

\

To prove the other inclusion observe that if for any $i$ holds:
$\phi_i \stackrel{G_r}{\sim}\psi_i$ then
 $(F_\bullet,\phi_\bullet)\stackrel{ G_r}{\sim}(F_\bullet,\psi_\bullet)$.
  Here $\prodl_i G_r$ acts on $\{F_i\}$ by a collection of automorphisms. Thus the statement follows directly from
part 1 of Theorem \ref{Thm.Results.Annihilators.for.GlGrGlr} and Theorem \ref{Thm.Background.Linearization}.

\

{\bf 2.} This now follows immediately from Lemma \ref{Thm.Proofs.Lemma.Chains.Morphisms}.
\epr

\subsection{Criteria for relative determinacy}\label{Sec.Proofs.Relative.Determinacy}
We start from a preparatory result.
\bel
Fix two $R$-modules $N\sseteq M$ and an ideal $J\sset R$.
\bee
\item $ann\quots{J\cdot M}{N}=ann(\quots{M}{N}):J$.
\item $ann\quots{M}{N}\cdot ann\quots{N}{J\cdot N}\sseteq ann\quots{M}{J\cdot N} \sseteq ann\quots{M}{N}\cap ann\quots{N}{J\cdot N}$.
\eee \eel \bpr 1. The proof is just the chain of immediate
equivalences: \beq \text{$\Big[f\in ann\quots{M}{N}:J\Big]$ iff
$\Big[f\cdot J\sseteq  ann\quots{M}{N}\Big]$ iff $\Big[f\cdot J\cdot
M\sseteq N\Big]$ iff $f\in ann\quots{J\cdot M}{N}$.} \eeq 2. If
$f\in ann\quots{M}{N}$ and $g\in ann\quots{N}{J\cdot N}$ then
$gf\cdot M\sseteq g\cdot N\sseteq J\cdot N$, i.e. $gf\in
ann\quots{M}{J\cdot N}$.

If $f\in ann\quots{M}{J\cdot N}$ then obviously $f\in
ann\quots{M}{N}$. As $N\sseteq M$ we get also: $f\cdot N\sseteq
J\cdot N$, i.e. $f\in ann\quots{N}{J\cdot N}$. \epr
\bex
In each bound of part 2 of the last lemma the equality is realizable.
 \bee[i.]
\item If $dim(R)$ is large enough and $J\sset R$ is generic for the given $M,N$ then
 $ann\quots{M}{N}\cdot  ann\quots{N}{J\cdot N}=ann\quots{M}{N}\cap ann\quots{N}{J\cdot N}$, hence part two is an equality.

\item Let $R=\k[[x,y]]$, with $J=\cm^p$ and $N=(x^p,y^p)\sset M=R$. Then $J\cdot N=\cm^{2p}$,
 $ann\quots{N}{J\cdot N}=\cm^p$,  $ann\quots{M}{N}=(x^p,y^p)$. Therefore:
\beq (x^p,y^p)= ann\quots{M}{N}\cap  ann\quots{N}{J\cdot
N}\supsetneq ann\quots{M}{J\cdot N}=\cm^{2p}= ann\quots{M}{N}\cdot
ann\quots{N}{J\cdot N}. \eeq \

\item Let  $J=N=(f)\sset M=R$. Then
$ann\quots{M}{N}\cdot ann\quots{N}{J\cdot N}=ann\quots{M}{J\cdot
N}=(f^2)\ssetneq ann\quots{M}{N}\cap ann\quots{N}{J\cdot N}=(f)$
\eee \eex
 This lemma implies the immediate

\bcor\label{Thm.Proofs.Corollary.For.Rel.Determinacy}
Let $\Si^{(J)}$, $G^{(I)}$ be as in \S\ref{Sec.Results.Relative.Determinacy}.
Suppose $T_{(\Si^{(J)},A)}=J\cdot T_{(\Si,A)}$ and $T_{(G^{(I)}A,A)}=I\cdot T_{(GA,A)}$.
Then:
\[
\Big(ann(T^1_{(\Si,G,A)})\cap ann\quotient{T_{(GA,A)}}{I\cdot
T_{(GA,A)}}\Big):J \supseteq
ann(T^1_{(\Si^{(J)},G^{(I)},A)})\supseteq\Big(ann(T^1_{(\Si,G,A)})\cdot
ann\quotient{T_{(GA,A)}}{I\cdot T_{(GA,A)}}\Big):J.
\]
\ecor
\bex
For $G=G_r$ and $I=R$ we get $ann(T^1_{(\Si^{(J)},G_r,A)})=ann.coker(A):J$.
If $R$ satisfies the condition \eqref{Eq.Intro.assumption.R.Noetherian.or.surjects}
 and $J^q\sseteq ann.coker(A)$ then $A+Mat_{m\times n}(J^{q+k})\sseteq G^{(J^k)}_rA$, for any $k>0$.
This strengthens and generalizes \cite[Theorem 5.2]{Cutkosky-Srinivasan}: ``For $R$ a domain and complete \wrt $I_m(A_1)$-filtration, there exist $k_1,k_2$ such that if $A_1\equiv A_2\ mod(I_m(A_1))^j$ for $j\ge k_1$ then $A_1=UA_2V^{-1}$, with $U=\one\ mod(I_m(A_1))^{j-k_2}$, $V=\one\ mod(I_m(A_1))^{j-k_2}$."
\eex
\bex
Let $A\in Mat^{sym}_{m\times m}(\cm)$ where $R$ satisfies condition \eqref{Eq.Intro.assumption.R.Noetherian.or.surjects}
and is  Henselian.
 If $J^q\sseteq I_m(A):I_{m-1}(A)$ then $\{A\}+Mat^{sym}_{m\times m}(J^{q+k})\sseteq G^{(J^k)}_{congr}A$ for any $k>0$.
\eex


\begin{thebibliography}{99}

\bibitem[Ahmed-Ruas]{Ahmed-Ruas} I.Ahmed, M.A. Soares-Ruas, {\em Determinacy of Determinantal Varieties},  arXiv:1710.11189.

\bibitem[Arnol'd68]{Arnold68} V.I. Arnol'd, {\em Singularities of smooth mappings.} Uspehi Mat. Nauk 23 1968 no. 1, 3--44.

\bibitem[AGLV-1]{AGLV} V.I. Arnol'd, V.V. Goryunov, O.V. Lyashko, V.A. Vasil'ev, {\em Singularity theory. I.}
Reprint of the original English edition from the series Encyclopaedia of Mathematical Sciences
[ Dynamical systems. VI, Encyclopaedia Math. Sci., 6, Springer, Berlin, 1993]. Springer-Verlag, Berlin, 1998. iv+245 pp.
ISBN: 3-540-63711-7



\bibitem[Beasley-Laffey]{Beasley-Laffey} L.B. Beasley, T.J. Laffey, {\em Linear operators on matrices: the invariance of rank-k matrices.} Linear Algebra Appl. 133 (1990), 175--184.


\bibitem[B.K.1]{Belitski-Kerner} G. Belitski, D. Kerner, {\em Group actions on filtered modules and finite determinacy.
Finding large submodules in the orbit by linearization},  C. R. Math. Rep. Acad. Sci. Canada Vol. 38 (4) 2016, pp. 113--155 .

\bibitem[B.K.2]{Belitski-Kerner2} G. Belitski, D. Kerner, {\em Finite determinacy of matrices over local rings.II.
 Tangent modules to the miniversal deformations for group-actions involving the ring automorphisms},  arXiv:1604.06247.

\bibitem[B.K.3]{Belitski-Kerner3} G. Belitski, D. Kerner, in preparation.




\bibitem[Birkhoff-MacLane]{Birkhoff-MacLane} S. Mac Lane, G. Birkhoff, {\em Algebra.} Third edition. Chelsea Publishing Co., New York, 1988. xx+626 pp.


\bibitem[Bruce]{Bruce2003} J.W. Bruce, {\em On families of symmetric matrices.} Dedicated to Vladimir I. Arnold
on the occasion of his 65th birthday. Mosc. Math. J. 3 (2003), no. 2, 335--360, 741


\bibitem[Bruce-Goryunov-Zakalyukin]{Bruce-Goryun-Zakal02} J.W. Bruce, V.V. Goryunov, V.M. Zakalyukin,
{\em Sectional singularities and geometry of families of planar quadratic forms}. Trends in singularities, 83--97,
Trends Math., Birkh\"{a}user, Basel, 2002


\bibitem[Bruce-Tari]{Bruce-Tari04} J.W. Bruce, F. Tari, {\em On families of square matrices}.
Proc. London Math. Soc. (3) 89 (2004), no. 3, 738--762.



\bibitem[Buchsbaum-Eisenbud]{Buchsbaum-Eisenbud}, D.A. Buchsbaum; D. Eisenbud, {\em  Algebra structures for finite free resolutions, and some structure theorems for ideals of codimension 3}. Amer. J. Math. 99 (1977), no. 3, 447--485



\bibitem[Cutkosky-Srinivasan]{Cutkosky-Srinivasan} S.D. Cutkosky, H. Srinivasan, {\em Equivalence and finite determinancy of mappings.} J. Algebra 188 (1997), no. 1, 16--57.


\bibitem[Damon84]{Damon1984} J. Damon, {\em The unfolding and determinacy theorems for subgroups of A and K.} Mem. Amer. Math. Soc. 50 (1984), no. 306, x+88 pp


\bibitem[Damon-Pike.1]{Damon.Pike.1} J. Damon,  B. Pike, {\em   Solvable Groups, Free Divisors and Nonisolated Matrix Singularities I: Towers of Free Divisors},    arXiv:1201.1577, to appear in Annales de l'Institut Fourier.


\bibitem[Damon-Pike.2]{Damon.Pike.2} J. Damon,  B. Pike, {\em Solvable groups, free divisors and nonisolated matrix singularities II: Vanishing topology.} Geom. Topol. 18 (2014), no. 2, 911--962



\bibitem[Dieudonn\'{e}]{Dieudonne} J. Dieudonn\'{e}, {\em Sur une g\'{e}n\'{e}ralisation du groupe
 orthogonal \`{a} quatre variables.} Arch. Math. 1, (1949), 282--287.




\bibitem[de Jong-van Straten]{de Jong-van Straten1990} T. de Jong, D. van Straten, {\em A deformation theory for nonisolated singularities.}
 Abh. Math. Sem. Univ. Hamburg 60 (1990), 177--208







\bibitem[du Plessis1982]{du Plessis1982} A. du Plessis, {\em On the genericity of topologically finitely-determined map-germs.} Topology 21 (1982), no. 2, 131--156.

\bibitem[du Plessis1983]{du Plessis1983} A. du Plessis,  {\em Genericity and smooth finite determinacy.}
  Singularities, Part 1 (Arcata, Calif., 1981),  295--312, Proc. Sympos. Pure Math., 40, Amer. Math. Soc., Providence, RI, 1983




\bibitem[Eisenbud]{Eisenbud} D. Eisenbud, {\em Commutative algebra. With a view toward algebraic geometry.}
 Graduate Texts in Mathematics, 150. Springer-Verlag, New York, 1995. xvi+785


\bibitem[Elkik]{Elkik73} R. Elkik, {\em Solutions d'\'{e}quations \`{a} coefficients dans un anneau hens\'{e}lien.} Ann. Sci. \'{E}cole Norm. Sup. (4) 6 (1973), 553-603 (1974)






\bibitem[Fr\"{u}hbis-Kr\"{u}ger.Zach]{Fr.Kr-Za.2015} A. Fr\"{u}hbis-Kr\"{u}ger, M. Zach, {\em On the Vanishing Topology of Isolated Cohen-Macaulay Codimension 2 Singularities}, arXiv: 1501.01915 .



\bibitem[Goryunov-Mond]{Goryun-Mond05} V. Goryunov, D. Mond, {\em Tjurina and Milnor numbers of matrix singularities.}
 J. London Math. Soc. (2) 72 (2005), no. 1, 205--224.

\bibitem[Goryunov-Zakalyukin]{Goryun-Zakal03} V.V. Goryunov, V.M. Zakalyukin,
{\em Simple symmetric matrix singularities and the subgroups of Weyl groups $A_\mu$, $D_\mu$, $E_\mu$}.
Dedicated to Vladimir I. Arnold on the occasion of his 65th birthday. Mosc. Math. J. 3 (2003), no. 2, 507--530, 743--744


\bibitem[Grandjean00]{Grandjean00} V. Grandjean, {\em Finite determinacy relative to closed and finitely generated ideals.} Manuscripta Math. 103 (2000), no. 3, 313--328.


\bibitem[Grandjean04]{Grandjean04} V. Grandjean, {\em Infinite relative determinacy of smooth function germs with transverse isolated singularities and relative ${\L}$ojasiewicz conditions.} J. London Math. Soc. (2) 69 (2004), no. 2, 518--530.





\bibitem[Greub]{Greub} W. Greub, {\em Multilinear algebra.} Second edition. Universitext. Springer-Verlag, New York-Heidelberg, 1978. vii+294 pp.

\bibitem[Greuel-Pham]{Greuel-Pham} G.-M. Greuel, T.H. Pham
 {\em Mather-Yau Theorem in Positive Characteristic}, arXiv:1308.5153.


\bibitem[Haslinger]{Haslinger} G.-J. Haslinger,  {\em Families of Skew-symmetric Matrices}
  Ph.D. Thesis, The University of Liverpool (United Kingdom). 2001.



\bibitem[Heymans]{Heymans} P. Heymans, {\em Pfaffians and skew-symmetric matrices.}
Proc. London Math. Soc. (3) 19 1969 730--768.


\bibitem[Huneke-Swanson]{Huneke-Swanson} C. Huneke, I. Swanson, {\em Integral closure of ideals, rings, and modules.}
  London Mathematical Society Lecture Note Series, 336. Cambridge University Press, Cambridge, 2006. xiv+431 pp.




\bibitem[Jiang]{Jiang} G. Jiang, {\em Functions with non-isolated singularities on Singular Spaces}, PhD thesis, Utrecht Netherlands, 1997.

\bibitem[Kucharz]{Kucharz} W. Kucharz, {\em Power series and smooth functions equivalent
 to a polynomial.} Proc. Amer. Math. Soc. 98 (1986), no. 3, 527--533.






\bibitem[Mather]{Mather1968} J.N. Mather,
{\em Stability of $C^\infty$ mappings. I. The division theorem.} Ann. of Math. (2) 87 1968 89--104.
\\{\em Stability of $C^\infty$ mappings. II. Infinitesimal stability implies stability.} Ann. of Math. (2) 89 1969 254--291.
\\{\em Stability of $C^\infty$ mappings. III. Finitely determined map-germs.} Publ. Inst. Hautes \'{E}tudes Sci. Publ. Math. No. 35, 1968, 279--308.
\\{\em Stability of $C^\infty$ mappings. IV. Classification of stable germs by R-algebras.} Publ. Inst. Hautes \'{E}tudes Sci. Publ. Math. No. 37 1969 223--248.
\\{\em Stability of $C^\infty$ mappings. V. Transversality.} Advances in Math. 4 1970 301--336 (1970).
\\ {\em Stability of $C^\infty$ mappings. VI: The nice dimensions.}
 Proceedings of Liverpool Singularities-Symposium, I (1969/70), pp. 207--253. Lecture Notes in Math., Vol. 192,
 Springer, Berlin, 1971




\bibitem[Li-Pierce]{Li-Pierce} C.-K. Li, S. Pierce, {\em Linear preserver problems.}
Amer. Math. Monthly 108 (2001), no. 7, 591--605.


\bibitem[Malgrange]{Malgrange} B. Malgrange, {\em Ideals of differentiable functions.} Tata Institute of Fundamental Research Studies in Mathematics, No. 3 Tata Institute of Fundamental Research, Bombay; Oxford University Press, London 1967 vii+106 pp



\bibitem[Moln\'{a}r]{Molnar} L. Moln\'{a}r, {\em Selected preserver problems on algebraic structures of linear
 operators and on function spaces.} Lecture Notes in Mathematics, 1895. Springer-Verlag, Berlin, 2007.



\bibitem[Pellikaan88]{Pellikaan88} R. Pellikaan, {\em
 Finite determinacy of functions with nonisolated singularities.}
 Proc. London Math. Soc. (3) 57 (1988), no. 2, 357--382.


\bibitem[Pellikaan90]{Pellikaan} R. Pellikaan {\em Deformations of hypersurfaces with a one-dimensional singular locus.} J. Pure Appl. Algebra 67 (1990), no. 1, 49--71.

\bibitem[Pereira.10]{Pereira.phd} M.S. Pereira, {\em Variedades Determinantais e Singularidades de Matrizes}, PhD thesis, 2010,
 Universidade de S\~{a}ão Paulo, USP, Brasil.

\bibitem[Pereira.17]{Pereira.2017} M.S. Pereira, {\em Properties of G-Equivalence of Matrices},  arXiv:1711.02156.

\bibitem[Pham]{Pham} T.H. Pham {\em On finite determinacy of hypersurface singularities and matrices in arbitrary characteristic}, PhD thesis, Technische Universit\"{a}t Kaiserslautern, 2016.






\bibitem[Rudin]{Rudin-book} W.Rudin, {\em Real and complex analysis.} Third edition.
McGraw-Hill Book Co., New York, 1987. xiv+416 pp

\bibitem[Siersma83]{Siersma83} D. Siersma, {\em Isolated line singularities. Singularities}, Part 2 (Arcata, Calif., 1981), 485--496,
Proc. Sympos. Pure Math., 40, Amer. Math. Soc., Providence, RI, 1983

\bibitem[Siersma00]{Siersma} D. Siersma, {\em The vanishing topology of non isolated singularities.} New developments
in singularity theory (Cambridge, 2000), 447--472, NATO Sci. Ser. II Math. Phys. Chem., 21, Kluwer Acad. Publ.,
Dordrecht, 2001





\bibitem[Teissier]{Teissier} B.Teissier, {\em Vari\'{e}t\'{e}s polaires. II. Multiplicit\'{e}s polaires, sections planes, et conditions de Whitney.}      Algebraic geometry (La R\'{a}bida, 1981), 314--491, Lecture Notes in Math., 961, Springer, Berlin, 1982.


\bibitem[Thilliez]{Thilliez}  V. Thilliez, {\em Infinite determinacy on a closed set for smooth germs with non-isolated singularities}, Proc. Amer. Math. Soc. 134 (2006), 1527--1536

\bibitem[Tougeron]{Tougeron1968} J.C. Tougeron, {\em Id\'{e}aux de fonctions diff\'{e}rentiables. I.}
 Ann. Inst. Fourier (Grenoble) 18 1968 fasc. 1, 177--240


\bibitem[Wall79]{Wall-1979} C.T.C. Wall, {\em Are maps finitely determined in general?}
Bull. London Math. Soc. 11 (1979), no. 2, 151--154.


\bibitem[Wall81]{Wall-1981} C.T.C. Wall, {\em Finite determinacy of smooth map-germs.} Bull. London
Math. Soc. 13 (1981), no. 6, 481--539.

\bibitem[Wilson81]{Wilson81} L.-C. Wilson, {\em Infinitely determined map germs.} Canad. J. Math. 33 (1981), no. 3, 671--684.

\bibitem[Wilson82]{Wilson82} L.-C. Wilson, {\em Map-germs infinitely determined with respect to right-left equivalence.}
Pacific J. Math. 102 (1982), no. 1, 235--245.
\end{thebibliography}
\end{document}